\documentclass[12pt,leqno,amssymb]{amsart}

\usepackage{amsmath,amssymb}

\usepackage{amsthm}
\usepackage[english,francais]{babel}

\numberwithin{equation}{section} 

\def\emty{\emptyset}
\def\cbar{\widehat{\C}}
\def\dis{\displaystyle}
\newcommand\ben{\begin{enumerate}}
\newcommand\een{\end{enumerate}}
\newcommand\bit{\begin{itemize}}
\newcommand\eit{\end{itemize}}
\def\lag{\langle}
\def\rag{\rangle}
\renewcommand\hat{\widehat}

\def\AAA{{\mathcal A}}
\def\BBB{{\mathcal B}}
\def\EEE{{\mathcal E}}
\def\HHH{{\mathcal H}}
\def\MMM{{\mathcal M}}
\def\NNN{{\mathcal N}}
\def\PPP{{\mathcal P}}
\def\SSS{{\mathcal S}}
\def\TTT{{\mathcal T}}

\def\supp{\mbox{\rm supp\,}}
\def\rond{\mbox{\rm Rond\,}}
\def\diam{\mbox{\rm diam\,}}
\def\dist{\mbox{\rm dist}}
\def\mod{\hbox{\rm mod}}
\def\th{th\'eor\`eme }
\def\Th{Th\'eor\`eme }
\def\qc{quasiconforme }
\def\homeo{hom\'eomorphisme }
\def\homeos{hom\'eomorphismes }
\def\qcs{quasiconformes }
\def\ccs{composantes connexes }

\def\Preuve{{\noindent\sc D\'emonstration.} }
\def\endp{\hspace*{\fill}$\rule{.55em}{.55em}$ \smallskip}

\def\al{\alpha}

\def\g{\gamma}
\def\G{\Gamma}
\def\si{\sigma}
\def\ep{\varepsilon}
\def\la{\lambda}
\def\De{\Delta}
\def\de{\delta}
\def\oo{\mbox{$\Omega$} }
\def\om{\mbox{$\omega$} }

\def\R{\mbox{$\mathbb R$}}
\def\C{\mbox{$\mathbb C$}}
\def\N{\mbox{$\mathbb N$}}
\def\SS{\mbox{$\mathbb S$}}
\def\HH{\mbox{$\mathbb H$}}
\def\DD{\mbox{$\mathbb D$}}
\def\PP{\mbox{$\mathbb P$}}

\def\hmod{\hat{\hbox{\rm mod}}}
\def\hrho{\hat{\rho}}
\def\hg{\hat{\gamma}}
\def\hG{\hat{\Gamma}}
\def\heta{\hat{\eta}}

\theoremstyle{plain}
\newtheorem{theo}{Th\'eor\`eme}[section]
\newtheorem{prop}[theo]{Proposition}
\newtheorem{coro}[theo]{Corollaire}
\newtheorem{conj}[theo]{Conjecture}
\newtheorem{defi}[theo]{D\'efinition}
\newtheorem{lemm}[theo]{Lemme}
\newtheorem{rema}[theo]{Remarque}

\newtheorem{ques}[theo]{Question}
\newtheorem{newthm}{Theorem}
\newcommand{\REFTHM}[1] { \begin{theo}\label{#1} }
\newcommand{\ENDTHM}{\end{theo}}
\newcommand{\REFPROP}[1]{\begin{prop}\label{#1} }
\newcommand{\ENDPROP}{\end{prop} }
\newcommand{\REFLEM}[1]{\begin{lemm}\label{#1} }
\newcommand{\ENDLEM}{\end{lemm} }
\newcommand{\REFCOR}[1]{\begin{coro}\label{#1} }
\newcommand{\ENDCOR}{\end{coro} }
\newcommand{\REFNTH}[1] { \begin{newthm}\label{#1} }
\newcommand{\ENDNTH}{\end{newthm}}

\date{Version du \today}

\title{Empilements de cercles et modules combinatoires}
\author{Peter Ha\"{\i}ssinsky}
\address{LATP/CMI\\ Universit\'e de Provence\\ 39, rue Fr\'ed\'eric Joliot-Curie\\  
13453 Marseille cedex 13\\France}
\email{phaissin@cmi.univ-mrs.fr}
\begin{document}

\maketitle

\noindent{\bf R\'esum\'e.} \hrulefill \\[.1cm]
Le but cette note est de tenter d'expliquer les liens \'etroits qui unissent
la th\'eorie des empilements de cercles et des modules combinatoires, et de comparer
les approches \`a la conjecture de J.W.\,Cannon qui en d\'ecoulent \cite{cs,cfp1,cfp2,bk0,bk1,bk2}.%

\noindent\hrulefill \\[.1cm]
\bigskip

\vspace*{-0.5cm}

\noindent{\bf Abstract.} \hrulefill \\[.1cm]
The aim of this article is to explain the deep relationships
between circle-packings and combinatorial moduli of curves, and to compare
the approaches to Cannon's conjecture to which they lead
 \cite{cs,cfp1,cfp2,bk0,bk1,bk2}.

\noindent\hrulefill \\[.1cm]

\tableofcontents

\newpage 

\section*{INTRODUCTION}

Plusieurs approches ont \'et\'e d\'evelopp\'ees pour estimer des modules de courbes dans des espaces
topologiques munis de structures tr\`es faibles. Les premiers r\'esultats dans ce sens sont certainement
d\^us \`a P.\,Pansu qui g\'en\'eralise la notion de module dans le cadre d'espaces munis de familles
(ind\'enombrables) de \og boules \fg\ \cite{pan}. D'autres m\'ethodes ont aussi \'et\'e d\'evelopp\'ees
par J.\,Heinonen et P.\,Koskela notamment pour \'etudier les \homeos \qcs et leurs avatars dans
des espaces de Carnot, voire dans des espaces m\'etriques plus g\'en\'eraux \cite{hk1}. 

L'objet de ce texte est d'unifier autant que possible les approches de J.W.\,Cannon, W.\,Floyd et W.\,Parry et
de M.\,Bonk et B.\,Kleiner \`a des probl\`emes d'uniformisation de surfaces motiv\'es par la conjecture
de Cannon qui affirme qu'un groupe hyperbolique dont le bord
est hom\'eomorphe \`a $\SS^2$ admet une action klein\'eenne cocompacte.
 Ces deux approches reposent
sur des notions de modules combinatoires l\'eg\`erement diff\'erentes qui conduisent \`a des r\'esultats
apparemment similaires. Le r\'esultat central de cet article est le \Th \ref{comp} qui fournit un socle
commun pour en d\'eduire d'une part une version faible du \th de Riemann combinatoire de J.W.\,Cannon
(Corollaire \ref{oricomp}) et d'autre part les r\'esultats de M.\,Bonk et B.\,Kleiner
(voir p.ex. le Corollaire \ref{th11.1}). Le point clef des d\'emonstrations que l'on pr\'esente
est d'utiliser, \`a l'instar de M.\,Bonk et B.\,Kleiner, les empilements de cercles pour traduire
naturellement les donn\'ees combinatoires en donn\'ees analytiques. 

Pour ce faire, plusieurs variantes de modules combinatoires pr\'esent\'es par ces
diff\'erents auteurs sont d\'evelopp\'ees dans ce texte, et on \'etablit quelques comparaisons 
entre elles. Nous adoptons le point de vue du lecteur familier des  notions et techniques
de th\'eorie g\'eom\'etrique des fonctions pour d\'emontrer les propri\'et\'es de ces modules.
Ces notions sont rappel\'ees bri\`evement dans l'Appendice \ref{tgf}, o\`u quelques d\'emonstrations
sont aussi esquiss\'ees pour insister sur le traitement commun de ces deux th\'eories.

Ce texte est organis\'e en quatre parties et trois appendices. Les deux premi\`eres parties fournissent 
des rappels concernant les empilements
de cercles d'une part, et les modules combinatoires de courbes d'autre part. La troisi\`eme partie
contient les r\'esultats proprement dits de l'article et la comparaison des deux m\'ethodes de discr\'etisation.
La derni\`ere partie contient une  discussion sommaire
sur les approches \`a la conjecture de Cannon. L'Appendice \ref{tgf} d\'ecrit succinctement les
notions de th\'eorie g\'eom\'etrique des fonctions qui sont utilis\'ees tout au long du texte, et 
qui servent de motivation. Dans l'Appendice \ref{ap:dis}, on montre que ces m\'ethodes de discr\'etisation
sont aussi int\'eressantes dans le cadre d'espaces m\'etriques plus g\'en\'eraux. Enfin, dans l'Appendice \ref{apphaus},
on \'etablit un crit\`ere de convergence pour des fonctions d\'efinies sur des sous-ensembles compacts d'un espace
compact m\'etrique qui converge au sens de Hausdorff.

\medskip

Les r\'esultats de ce texte sont essentiellement connus, 
mais il est esp\'er\'e qu'ils sont mis en perspective de mani\`ere profitable, et que les arguments pr\'esent\'es
seront trouv\'es plus simples. 

\medskip

{\noindent\sc Remerciements.} Je remercie J.\,Los, L.\,Reeves et H.\,Short pour les nombreuses discussions que nous avons eues
autour des th\`emes abord\'es dans ce texte, sans lesquelles ce travail n'aurait pu voir le jour. Je remercie aussi Hadrien
Lar\^ome qui a relu avec beaucoup d'attention une grande partie de ce manuscrit. Je suis reconnaissant au rapporteur pour ses
nombreux conseils qui m'ont permis d'am\'eliorer cet article.

\medskip

\noindent{\bf Notations.} L'ensemble des r\'eels positifs ou nul est d\'esign\'e par $\R_+$. 
Dans tout le texte, si $a$ et $b$ sont strictement positifs, on \'ecrit $a\asymp b$ s'il existe 
une constante universelle $u>1$ telle que $(1/u)\cdot a\le b\le u \cdot a$, $a\lesssim b$ ou $b\gtrsim a$ s'il existe une
constante universelle $C>1$ telle que $a\le C\cdot b$.

\section{Empilements de cercles}

Un empilement de cercles $\EEE$  de la sph\`ere de Riemann $\cbar$ est 
une collection de disques ferm\'es d'int\'erieurs deux \`a deux disjoints telle que les \ccs du compl\'ementaire
de la r\'eunion de ces disques soient bord\'ees par exactement trois arcs de cercle. 

A un empilement $\EEE$, on associe son graphe d'incidence $G(\EEE)$ dont les sommets sont les
disques et les ar\^etes sont les paires de disques dont les fermetures s'intersectent. Lorsque 
$\EEE$ est fini, il s'agit d'une triangulation que l'on peut voir plong\'ee dans la surface.

Nous \'enon\c{c}ons d\`es maintenant un \th de repr\'esentation conforme pour les empilements qui correspond
\`a une r\'eciproque du fait pr\'ec\'edent.

\begin{theo}[P.\,Koebe \cite{koe}]\label{andreev} Soit $\TTT$ une triangulation de la sph\`ere $\SS^2$. Il existe un unique
empilement de cercles $\EEE$ de $\cbar$ de graphe d'incidence isomorphe \`a $\TTT$, \`a transformation
homographique pr\`es. \end{theo}

Une d\'emonstration plus moderne figure dans \cite{cv1}. 

W.P.\,Thurston a propos\'e de d\'emontrer le \th de repr\'esentation conforme de Riemann en utilisant les
empilements de cercles. Soient $\oo\subset\C$ un domaine simplement connexe born\'e du plan, et $z_1,z_2\in\oo$
deux points distincts. Pour $r>0$, on note $\TTT_r$ la triangulation de $\C$ par triangles \'equilat\'eraux engendr\'ee
par $[0,r]\times\{0\}$. 
On consid\`ere la r\'eunion des triangles pleins enti\`erement contenus dans $\oo$, et on d\'esigne 
par $\oo_r$ la composante connexe de son int\'erieur qui contient $z_1$. Si $r$ est assez petit, alors
$z_2$ en fait aussi partie. Soit $\TTT_r(\oo)=\TTT_r\cap\overline{\oo}_r$ la restriction de la triangulation
$\TTT_r$ \`a $\overline{\oo_r}$. 
 On ajoute \`a $\TTT_r(\oo)$ un sommet (abstrait) et des ar\^etes qui le relient aux
sommets du bord de $\TTT_r(\oo)$. On obtient ainsi une triangulation de la sph\`ere, pour laquelle on peut appliquer
le \Th \ref{andreev} afin d'obtenir un empilement de cercles $\EEE_r$ dont le compl\'ementaire du disque unit\'e
$\C\setminus\DD$ est associ\'e au
sommet rajout\'e. On note $f_r:\TTT_r(\oo)\to\DD$ l'application simpliciale induite par les triangulations sous-jacentes.
Quitte \`a composer par un automorphisme de $\DD$, on peut supposer que $f_r(z_1)=0$ et $f_r(z_2)>0$. 
Le but est alors de montrer que $f_r$ converge vers l'application de Riemann quand $r$ tend vers $0$, ce qu'ont fait B.\,Rodin et D.\,Sullivan
\cite{rs}.

\begin{theo} Le chemin $(f_r)_r$ converge uniform\'ement sur les compacts de $\oo$ vers la transformation conforme
$f:\oo\to\DD$ telle que $f(z_1)=0$ et $f(z_2)>0$ lorsque $r$ tend vers $0$.\end{theo} 

Un point crucial de leur d\'emonstration est de montrer que la suite est uniform\'ement quasiconforme (voir Appendice \ref{tgf}). Cette 
propri\'et\'e est une cons\'equence du lemme suivant.

\begin{lemm}[du collier] Il existe une constante $C=C(N)>0$ telle que si un disque $D$
de rayon $r$ 
d'un empilement de cercles a $N$ voisins, alors le rayon $r'$ de chaque voisin v\'erifie $r'>C(N) r$.\end{lemm}

Une d\'emonstration figure dans \cite{rs}.

\section{Modules combinatoires}

La notion classique de modules de courbes est expos\'ee dans l'Appendice \ref{tgf}. On s'en inspire pour en d\'evelopper une th\'eorie combinatoire.
\subsection{Modules dans un graphe}

Soit $G$ un graphe connexe de sommets $\SSS$ et d'ar\^etes $\AAA$. On munit ce graphe d'une m\'etrique
qui rend chaque ar\^ete isom\'etrique au segment $[0,1]$. Une {\sl courbe} sur le graphe 
sera un sous-graphe connexe de $G$, donc ses sommets sont dans $\SSS$, et ses ar\^etes dans $\AAA$.

Pour $Q\ge  1$, on note $\MMM_Q(G)$ l'ensemble des 
applications $\rho:\SSS\to \R_+$ telles que $0<\sum\rho(s)^Q<\infty$ que l'on appelle {\sl poids} ou {\sl m\'etriques admissibles},
ou {\sl $Q$-admissibles} pour \^etre plus pr\'ecis.

Soit $K$ un ensemble de  sommets de $G$\,; la {\sl $\rho$-longueur} de $K$ est par d\'efinition 
$$\ell_\rho (K)=\sum_{s\in K} \rho(s)$$
et son $Q$-volume $$V_{Q,\rho}(K)=\sum_{s\in K} \rho(s)^Q\,.$$
Lorsque $K=\SSS$, on notera plus simplement $V_{Q,\rho}(\SSS)=V_{Q}(\rho)$.

Si $\G$ est une famille de courbes dans $G$ et $\rho$ est une m\'etrique $Q$-admissible, on d\'efinit 
$$L_\rho(\G,\SSS)=\inf_{\g\in\G} \ell_\rho(\g),$$
et son $Q$-module combinatoire par
$$\mod_Q(\G,G) = \inf_{\rho\in\MMM_Q(G)} \frac{V_{Q}(\rho)}{L_\rho(\G,\SSS)^Q}= \inf_{\rho\in\MMM_Q(G)} \mod_Q(\G,\rho,\SSS).$$
Si, pour une m\'etrique $\rho$, on a $L_\rho(\G,\SSS)=0$, alors on pose $\mod_Q(\G,\rho,\SSS)=+\infty$.

On peut remarquer que si $\la >0$, alors  $\mod_Q(\G,\la\rho,\SSS)=\mod_Q(\G,\rho,\SSS)$, ce qui permet de normaliser
les m\'etriques par des conditions du type $L_\rho(\G,\SSS)\ge 1$, ou encore $V_{Q,\rho}(G)=1$.

On a les propri\'et\'es suivantes\,:
\ben 
\item  si $\Gamma_{1} \subset \Gamma_{2}$, $\mod_{Q}( \Gamma_{1},G) \leq
\mod_{Q} (\Gamma_{2},G)$;
\item $\mod_{Q} ( \bigcup_{i=1}^{\infty} (\Gamma_{i},G) \leq
\sum_{i=1}^{\infty} \mod_{Q} (\Gamma_{i},G)$\,;
\item si $\Gamma_{1}$ et $\Gamma_{2}$ sont deux familles de courbes
telles que toute courbe $\gamma_{1}$ dans $\Gamma_{1}$  poss\`ede une
sous-courbe $\gamma_{2} \in \Gamma_{2}$, alors $\mod_{Q}(\Gamma_{1},G)
 \leq \mod_{Q} (\Gamma_{2},G)$.
\end{enumerate}

En supposant $\SSS$ fini et en normalisant les \'el\'ements de $\MMM_Q(G)$ pour que $\sum\rho(s)^Q=1$, il n'est pas
difficile de montrer l'existence d'une m\'etrique extr\'emale {\it i.e.}, pour laquelle l'infimum est
atteint. Voici une adaptation au cas discret d'un crit\`ere de A.\,Beurling d'extr\'emalit\'e (voir l'Appendice \ref{tgf}).

\begin{prop}[crit\`ere de Beurling]\label{beurling} Soient $G$ un graphe fini et connexe, 
$\G$ une famille de courbes et $Q>1$.
Une m\'etrique $\rho$ est extr\'emale si et seulement si il existe une sous-famille
$\G_0\subset\G$ telle que \ben
\item pour toute $\g\in\G_0$, $\ell_\rho(\g) = L_\rho(\G,\SSS)$\,;
\item pour toute $h:\SSS\to\R$ telle que, pour $\g\in\G_0$, $\sum_{s\in \SSS(\g)} h(s)\ge 0$,\\
on ait $\sum_{s\in\SSS} h(s)\rho(s)^{Q-1}\ge 0$, o\`u $\SSS(\g)$ d\'esigne les sommets $s\in\SSS$ qui intersectent $\g$.\een
De plus, cette m\'etrique est unique \`a normalisation pr\`es.\end{prop}

Il est pratique d'introduire \`a ce niveau quelques notations et de noter quelques observations.
Si $G=(\SSS,\AAA)$ est un graphe fini et $\G$ une famille de courbes, on notera  $\G(s)$  l'ensemble
des courbes de $\G$ qui contiennent l'\'el\'ement $s\in\SSS$\,; de m\^eme, $\SSS(\g)$ d\'esignera 
l'ensemble des sommets de $\SSS$ contenus dans la courbe $\g\in\G$. 
Remarquons que si $\rho$ est extr\'emale, alors pour tout sommet  $s$ du
support $\supp \rho$ de $\rho$, il existe une courbe $\g\in\G$
telle que $\ell_\rho(\g)=L_\rho(\G,\SSS)$, car sinon on pourrait diminuer la valeur de $\rho(s)$, ce qui diminuerait
$\mod_Q(\G,\rho,\SSS)$. 

Enfin, on consid\`ere l'espace vectoriel r\'eel $\R^{\SSS}$ muni de sa base canonique $\{e_s\}_{s\in\SSS}$
et du produit scalaire $\lag\cdot,\cdot\rag$ qui rend cette base orthonorm\'ee. A chaque courbe $\g\in\G$,
on associe le vecteur $u_{\g}=\sum_{s\in\SSS(\g)} e_s$. Une m\'etrique $\si$ est donc un vecteur de $\R^{\SSS}$
et si $\g\in\G$, alors $\ell_{\si}(\g)=\lag u_{\g},\si\rag$. On d\'efinit aussi 
$v(\si)=\sum_{s\in\SSS} |\si(s)|^Q$ et on note que $v:\R^{\SSS}\to\R$ est strictement convexe si $Q>1$ et  lin\'eaire si $Q=1$. 
En tant que fonction num\'erique d\'efinie sur l'espace vectoriel $\R^{\SSS}$, $v$ est diff\'erentiable et son
gradient $\nabla_\rho v$ en un vecteur $\rho$ a pour coordonn\'ees $Q \rho(s)|\rho(s)|^{Q-2}$, $s\in\SSS$.

\begin{lemm}\label{lbeur} Avec les notations ci-dessus et $Q>1$, si $\rho$ est extr\'emale, alors il existe une sous-famille
$\G_0\subset\G$ et un vecteur $(\la_{\g})\in (\R_+)^{\G_0}$ tels que \ben
\item pour toute $\g\in\G_0$, $\ell_\rho(\g) = L_\rho(\G,\SSS)$\,;
\item $\nabla_{\rho}v = \sum_{\g\in\G_0} \la_\g u_\g$.\een
\end{lemm}

Ce lemme appara\^{\i}t dans \cite{cfp1} dans le cas $Q=2$.

\medskip

\Preuve Tout d'abord, on note $\G_0$ l'ensemble des courbes de $\G$ de longueur minimale, que l'on fixe
\'egale \`a $1$. Notons $\de=\inf_{\G\setminus\G_0}\ell_\rho(\g)-1$. Puisque l'on a un nombre fini de
vecteurs $u_\g$ possibles, on a $\de>0$. 

On consid\`ere l'ensemble $\De$ des vecteurs $\si\in\R^{\SSS}$ tels que $\ell_\si(\g)=1$ pour toutes les courbes
$\g\in\G_0$. Il s'agit
d'un espace affine, donc la restriction de $v$ \`a $\De$ admet un unique minimum en un point $\rho_1$ car $v$
est strictement convexe.
Montrons que $\rho_1=\rho$ par l'absurde. Pour $t\in[0,1]$, on note $\rho_t=\rho+ t(\rho_1-\rho)$.
Donc, si $\g\in\G_0$ alors $\ell_{\rho_t}(\g)=1$, et sinon 
$$\ell_{\rho_t}(\g)\ge (1-t)(1+\de) + t\ell_{\rho_1}(\g).$$ Du coup, si $t$ est assez petit, alors
on a $L_{\rho_t}(\G,\SSS)\ge 1$, et par convexit\'e, on devrait avoir $v(\rho_t) <v(\rho)$, ce qui
contredit l'extr\'emalit\'e de $\rho$. Par suite, $\rho$ est la solution du probl\`eme variationnel qui consiste
\`a minimiser $v$ sur $\De$.

La th\'eorie des extrema li\'es indique que $\nabla_{\rho} v$ est orthogonal \`a $\De$, autrement dit,
il existe des multiplicateurs de Lagrange $(\la_\g)\in\R^{\G_0}$ tels que 
$$\nabla_{\rho} v=\sum_{\g\in\G_0}\la_\g u_\g.$$

Pour conclure, il reste \`a montrer que l'on peut choisir ces multiplicateurs positifs ou nuls.
On proc\`ede par l'absurde  en montrant que $\rho$ n'est pas extr\'emale s'il existe
un multiplicateur strictement n\'egatif.
On note $E$ le sous-espace de $\R^{\SSS}$ engendr\'e par $\{u_\g,\ \g\in\G_0\}$ et  $F$
l'orthogonal de $\nabla_{\rho} v$ dans $E$ et on d\'efinit l'espace affine $A=\nabla_\rho v+F$.

L'espace $A$ s'identifie avec un sous-espace de l'espace projectif $\PP(E)$ de $E$. Puisque le produit
scalaire des $u_\g$ par $\nabla_\rho v$ est toujours positif, chaque demi-droite de vecteur directeur
$u_\g$ coupe $F$ en un point $p_\g$, et on note $0$ le point qui correspond \`a $\nabla_\rho v$.

Si $\nabla_{\rho} v$ n'est pas dans le c\^one engendr\'e par $\{u_\g\}_{\g\in\G_0}$,
alors $0$ n'est pas dans l'enveloppe convexe des points $p_\g$. Du coup, ces points sont
dans un demi-espace de $A$ et 
il existe un vecteur $\si\in F$ dont le produit scalaire avec tous les $u_\g$ est strictement positif, et m\^eme minor\'e 
par une constante $\eta>0$ , puisqu'ils sont en nombre fini. 

Posons $\rho_t=\rho + t\si$ avec $t>0$. Pour $\g\in\G_0$, on a 
$$\ell_{\rho_t}(\g)=\ell_\rho(\g) + t\lag u_\g, \si\rag\ge 1+ t\eta\,.$$
Et si $\g\in\G\setminus\G_0$, $\ell_{\rho_t}(\g) \ge 1+\de  +t\lag u_\g, \si\rag \ge 1+ t\eta$ si $t$ est assez petit,
donc $L_{\rho_t}(\G)\ge 1 +\eta t$  si $t$ est assez petit. 

 Si $t$ assez petit, alors pour tout $s\in\SSS$ dans
le support de $\rho$, on a $\rho_t(s)>0$. Si $\rho(s)=0$, alors $|\rho_t(s)|^Q= o(t)$. Par cons\'equent
$$\begin{array}{ll}v(\rho_t)& =\dis\sum _{s\in\SSS} |\rho(s) + t\si(s)|^Q\\ &\\
& =\dis\sum _{s\in\hbox{\small supp}\,\rho}\rho(s)^Q ( 1+ Q t\si(s) \rho(s)^{-1}+o(t)) + o(t)\\ &\\
& = v(\rho) +\dis\sum _{s\in\SSS}  Q \rho(s)^{Q-1}t\si(s)+o(t)\\ &\\
& = v(\rho) + t\lag\nabla_{\rho} v, \si\rag +o(t)\\ &\\
& =v(\rho) +o(t)\,,\end{array}$$
car $\nabla_{\rho} v$ et $\si$ sont orthogonaux par d\'efinition.

Donc $$\mod_Q(\G,\rho_t,\SSS) \le\frac{v(\rho)+o(t)}{(1+\eta t)^Q}= \mod_Q(\G,\rho,\SSS)(1-Q\eta t+o(t))+o(t)$$
ce qui contredit l'extr\'emalit\'e de $\rho$. Donc $\nabla_{\rho} v$ est dans le c\^one engendr\'e par 
$\{u_\g\}_{\g\in\G_0}$ et il existe $(\la_{\g})\in (\R_+)^{\G}$ tel que  
$\nabla_{\rho} v= \sum_{\g\in\G_0} \la_\g u_\g$.\endp

\medskip

{\noindent\sc D\'emonstration de la Proposition \ref{beurling}.}
Supposons que $\rho$ v\'erifie les conditions de A.\,Beurling. Soit $\si\in\MMM_Q(G)$ que l'on suppose normalis\'ee pour que $L_{\si}(\G)=L_\rho(\G)$. Du coup, 
$\ell_\si(\g)\ge L_\rho(\G)=\ell_\rho(\g)$ pour toute $\g\in\G_0$. Donc, en posant $h=\si-\rho$, on obtient
$$\sum_{s\in\SSS} (\si-\rho) \rho^{Q-1}\ge 0$$
soit $$\sum_{s\in\SSS} \rho^Q\le \sum_{s\in\SSS}\si\rho^{Q-1} \le \left(\sum_{s\in\SSS}\si^Q\right)^{1/Q}\cdot\left(\sum_{s\in\SSS} \rho^Q\right)^{1-1/Q}$$
par l'in\'egalit\'e de H\"older. Du coup, on a $\mod_Q(\G,\rho,\SSS)\le\mod_Q(\G,\si,\SSS)$.
Le cas d'\'egalit\'e se produit lorsque l'in\'egalit\'e de H\"older est une \'egalit\'e, donc si $\si^Q$ et 
$(\rho^{Q-1})^{Q/(Q-1)}$
sont proportionnelles, soit si $\si=\rho$ d'apr\`es leurs normalisations.

\bigskip

R\'eciproquement, supposons que $\rho$ est extr\'emale et normalisons-la pour que $L_\rho(\G)=1$.
D'apr\`es  le lemme pr\'ec\'edent, il existe $\G_0\subset\G$ et $(\la_{\g})\in (\R_+)^{\G}$ tels que 
\ben
\item pour toute $\g\in\G_0$, $\ell_\rho(\g) = L_\rho(\G,\SSS)$\,;
\item $\nabla_{\rho} v= \sum_{\g\in\G_0} \la_\g u_\g$.\een 

Si $h:\SSS\to\R$ v\'erifie $\sum_{s\in \SSS(\g)} h(s)\ge 0$ pour toute courbe $\g\in\G_0$, alors
$$\begin{array}{ll} \dis\sum_{s\in\SSS} h(s)\rho(s)^{Q-1} & =\dis\sum_{s\in\SSS} h(s)\left(\dis\frac{1}{Q}\dis\sum_{\g\in\G_0(s)} \la_\g\right)  \\ &\\
& = \dis\frac{1}{Q}\dis\sum_{\g\in\G_0} \la_\g\left(\dis\sum_{s\in\SSS(\g)} h(s)\right)\\ & \\
& \ge 0\,.\end{array}$$
Ceci \'etablit la r\'eciproque.
\endp

{\noindent\bf Condensateurs dans un graphe.} Si $E$ et $F$ sont deux ensembles
de sommets disjoints d'un graphe $G$, on consid\`ere la famille $\G(E,F,G)=\G(E,F)$
des courbes qui joignent $E$ \`a $F$ et on d\'efinit le module,
aussi appel\'e capacit\'e, par
$$\mod_Q(E,F,G)=\mod_Q(\G(E,F),G)\,.$$
On dira que $(E,F)$  est un condensateur si $E$ et $F$ sont les sommets
de sous-graphes connexes de $G$.

La capacit\'e d\'epend de mani\`ere essentielle
des ar\^etes du graphe.
Cependant, on a le r\'esultat suivant.

\begin{prop}\label{mgr} Soient  $\SSS$ un ensemble de sommets et $\AAA_1$ et $\AAA_2$ 
deux familles d'ar\^etes telles que les graphes associ\'es $G_1=(\SSS,\AAA_1)$ 
et $G_2=(\SSS,\AAA_2)$ soient connexes. Pour $j=1,2$, on munit $G_j$ de la distance de longueur
$|\cdot|_j$ pour laquelle chaque ar\^ete est de longueur $1$. On suppose que $\AAA_2\subset\AAA_1$,
$G_2$ est de valence born\'ee par $K$, et  que si
$|v_1-v_2|_{1}= 1$ alors $|v_1-v_2|_{2}\le k$. Soit $Q\ge 1$.

{\rm (1)} Pour toute paire de sommets
disjoints $(E,F)$,
on a $$\mod_Q(E,F,G_1)\asymp \mod_Q(E,F,G_2)\,,$$ 
o\`u les constantes implicites ne d\'ependent que de $k$, $K$ et $Q$.

{\rm (2)} Si $(E_1,F_1)$ est un condensateur de $G_1$, tel que $E_1$ et $F_1$ sont \`a distance au moins $2k+1$ 
dans $G_1$, alors
il existe un condensateur $(E_2,F_2)$ de $G_2$, tel que $E_1\subset E_2$, $F_1\subset F_2$
et $$\mod_Q(E_1,F_1,G_1)\asymp \mod_Q(E_2,F_2,G_2)\,,$$ 
o\`u les constantes implicites ne d\'ependent que de $k$, $K$ et $Q$.

\end{prop}

\Preuve (1) 
On note $\G_1=\G(E,F,G_1)$ la famille de courbes associ\'ee \`a $\{E,F\}$ dans $G_1$ et $\G_2=\G(E,F,G_2)$ dans $G_2$. 

Tout d'abord, on a $\mod_Q(\G_2,G_2)\le \mod_Q(\G_1,G_1)$ car il y a plus de chemins dans $G_1$. 

R\'eciproquement, on se donne une m\'etrique $\rho_2:\SSS\to\R_+$ telle que $L_{\rho_2}(\G_2,G_2)=1$. 
On d\'efinit $\rho_1:\SSS\to\R_+$
par $\rho_1(s)=\max\{\rho_2(s'),\ s'\in B_2(s,k)\}$, o\`u $B_2(s,k)$ d\'esigne la boule centr\'ee en $s$ de rayon
$k$ pour la distance $|\cdot|_2$. Notons $N$ le cardinal maximal d'une boule combinatoire $B_2(s,k)$.
Si $\g_1$ est une courbe de $\G_1$, on lui associe une courbe $\g_2$ de
$\G_2$ passant par au moins
les m\^emes sommets, et dont les d\'etours sont de longueur au plus $N$. Du coup, pour chaque sommet $s$ de $\g_1$, on
a $\rho_1(s)\ge (1/N)\ell_{\rho_2}(\g_2\cap B_2(s,k))$ et on en d\'eduit que $L_{\rho_1}(\G_1,G_1)\ge 1/N$.

De m\^eme, on a $$V_Q(\rho_1)\le\sum_{s\in \SSS} \sum_{v\in B_2(s,k)}\rho_2(v)^Q\le N \sum_{s\in \SSS} \rho_2(s)^Q\le N\cdot V_Q(\rho_2)$$
car chaque $v$ ne peut appara\^{\i}tre qu'au plus $N$ fois. On en d\'eduit $\mod_Q(\G_2,G_2)\gtrsim \mod_Q(\G_1,G_1)$.

(2) Si $E_1$ n'est pas connexe dans $G_2$, alors il existe deux sommets 
$s,s'\in E_1$ tels que $|s-s'|_1=1$ mais
$|s-s'|_2 >1$. On construit une courbe dans $G_2$ de longueur au plus $k$ qui relie $s$ \`a $s'$
que l'on ajoute \`a $E_1$. On continue ainsi jusqu'\`a obtenir $E_2\supset E_1$ comme sommets
d'un sous-graphe connexe de $G_2$ dans le $k$-voisinage de $E_1$. On proc\`ede de m\^eme avec
$F_1$ pour obtenir $F_2$. Par hypoth\`eses, $E_2$ et $F_2$ sont disjoints donc $(E_2,F_2)$ est un
condensateur de $G_2$. D'apr\`es (1), il suffit de montrer que $\mod_Q(E_1,F_1,G_2)\asymp \mod_Q(E_2,F_2,G_2)$.

Par construction, chaque courbe de $\G_1=\G(E_1,F_1,G_2)$ contient une courbe de $\G_2=\G(E_2,F_2,G_2)$, donc
$\mod_Q(E_1,F_1,G_2)\le \mod_Q(E_2,F_2,G_2)$. Pour la r\'eciproque, on consid\`ere  
une m\'etrique $\rho_1:\SSS\to\R_+$ telle que $L_{\rho_1}(\G_1,G_2)=1$. 
On d\'efinit $\rho_2:\SSS\to\R_+$
par $\rho_2(s)=\max\{\rho_1(s'),\ s'\in B_2(s,k)\}$.

Si $\g_2$ est une courbe de $\G_2$ d'extr\'emit\'es $e_2\in E_2$ et $f_2\in F_2$, on construit
des courbes $\g_e$ et $\g_f$ de longueur combinatoire au plus $k$ qui relient $e_2$ \`a $E_1$ 
et $f_2$ \`a $F_1$ respectivement. On concat\`ene alors $\g_2$ \`a $\g_e$ et $\g_f$ pour obtenir
une courbe $\g_1$ de $\G_1$. 

On remarque que
$$\rho_2(e_2)\ge  (1/N)\sum_{s\in B_2(e_2,k)}\rho_1(s)\ge (1/N)\ell_{\rho_1}(\g_e)$$
et $\rho_2(f_2)\ge (1/N)\ell_{\rho_1}(\g_f)$. Du coup, 
$\ell_{\rho_2}(\g_2)\ge (1/N) \ell_{\rho_1}(\g_1)\ge 1/N$. On conclut comme au point (1).\endp

\begin{defi}[graphes \'equivalents]\label{geq} Sous les conditions de la Proposition \ref{mgr}, on 
dira que $G_1$ contient le graphe \'equivalent $G_2$.\end{defi}

\subsection{Modules sur une surface}

Soit $\SSS$ un recouvrement localement fini d'un espace m\'etrique $X$. 
Ses \'el\'ements seront appel\'es  parfois des pi\`eces, et la maille de $\SSS$ sera le plus grand 
diam\`etre de ses pi\`eces.
Le {\sl nerf} $\NNN(\SSS)$ de  $\SSS$ est le graphe dont les sommets sont les
\'el\'ements de $\SSS$ et les ar\^etes sont constitu\'ees des paires de pi\`eces
qui s'intersectent. On dit qu'un recouvrement est {\sl de valence born\'ee} si
son nerf est de valence born\'ee, et {\sl l'indice de recouvrement} entre deux recouvrements
est le nombre maximal de pi\`eces de l'un des recouvrements qui intersectent
une m\^eme pi\`ece de l'autre. {\sl L'indice d'auto-recouvrement} est l'indice 
de recouvrement d'un recouvrement par rapport \`a lui-m\^eme. En particulier,
un graphe de valence born\'ee  a un indice d'auto-recouvrement fini.
On dira dans la suite qu'une suite de recouvrements est {\sl de valence born\'ee}
si la borne ne d\'epend pas du recouvrement de la suite.
Si $K\subset X$, on note $\SSS(K)$ l'ensemble des $s\in\SSS$ tels que $s\cap K\ne \emty$.

Par pavage, on entendra un recouvrement localement fini par compacts connexes d'int\'erieurs deux \`a deux disjoints.

Une courbe dans un espace topologique $X$ est une application continue $\g:[0,1]\to X$. On ne
distinguera pas en g\'en\'eral l'application $\g$ en tant que telle de son image dans $X$.

Si $\G$ est une famille de courbes dans $X$, on lui associe une famille de courbes $\widehat{\G}$ dans 
$\NNN(\SSS)$ comme suit: chaque courbe $\g$ traverse un ensemble de  pi\`eces adjacentes $\SSS(\g)$;
on associe \`a $\g$ cet ensemble $\SSS(\g)$, ainsi que les ar\^etes qui les connectent dans le graphe\,: 
cela d\'efinit une courbe $\hg$ de $\NNN(\SSS)$ au sens du paragraphe pr\'ec\'edent.
On note $\hG$ l'ensemble des courbes ainsi obtenues, qui correspond \`a $\{\SSS(\g),\g\in\G\}$ aux ar\^etes pr\`es.
On d\'efinit 
$$\mod_Q(\G,\SSS)=\mod_Q(\widehat{\G},\NNN(\SSS)).$$

Si $(E,F)$ est un condensateur {\it i.e.}, $E$ et $F$ sont deux continua
(compacts connexes non d\'eg\'en\'er\'es) disjoints, on d\'efinit deux modules. 
Le premier est associ\'e
\`a la famille des courbes $\G(E,F)$ qui joignent $E$ \`a $F$ dans $X$, 
et le second au condensateur vu dans $\NNN(S)$\,:
on consid\`ere $\widehat{E}=\SSS(E)$ et $\widehat{F}=\SSS(F)$ 
comme sous-ensemble de $\NNN(\SSS)$.

On d\'efinit alors $$\left\{\begin{array}{l}
\mod_Q(E,F,\SSS)=\mod_Q(\G(E,F),\SSS)\\
\hmod_Q(E,F,\SSS)=\mod_Q(\widehat{E},\widehat{F},\NNN(\SSS))\end{array}\right.$$

En g\'en\'eral, l'in\'egalit\'e suivante est satisfaite\,:
 $$
\mod_Q(E,F,\SSS)\le
\hmod_Q(E,F,\SSS).$$
\medskip

\begin{rema} Le module le plus maniable est bien s\^ur $\hmod_Q$ puisqu'il ne d\'epend que du nerf du recouvrement.
En revanche, il peut ne pas porter d'information si par exemple on n'a pas de courbes
entre les continua\,! Du coup, la signification de ce module se fera \`a travers l'autre, qui lui, d\'epend plus
pr\'ecis\'ement de la topologie de $X$.\end{rema}

\medskip

On suppose dor\'enavant que $X$ est une surface topologique.

Rappelons qu'un quadrilat\`ere est la donn\'ee d'un domaine de Jordan avec
quatre points distincts sur le bord d\'efinissant ses quatre c\^ot\'es.
Voici l'analogue du \Th \ref{andreev} ({\it cf.} \cite{cfp1,sch}).

\begin{theo} Soit $Q$ un quadrilat\`ere recouvert par un pavage $\SSS$ dont le nerf forme une triangulation
et soit 
$\G$ la famille de courbes qui joint deux c\^ot\'es oppos\'es donn\'es. Il existe
un rectangle $R$ et un pavage  de nerf identique au pr\'ec\'edent par carr\'es qui respecte les sommets. La longueur des ar\^etes
des carr\'es produit la m\'etrique extr\'emale pour $\mod_2(\G,\SSS)$. \end{theo}

\begin{rema} Si le pavage ne forme pas une triangulation, alors le pavage par carr\'es peut avoir un nerf sensiblement diff\'erent 
de l'original si plusieurs pi\`eces de $\SSS$ s'intersectent en un m\^eme point.\end{rema}

Sur une surface, il est pratique de travailler avec des anneaux, c'est-\`a dire des
 surfaces hom\'eomorphes
\`a $]0,1[\times \SS^1$. On leur associe les modules suivants.

\begin{defi}[modules combinatoires d'un anneau]\label{defi:modcombanneau} 
Si $A$ un anneau d'une surface recouverte par $\SSS$, 
on consid\`ere $\G_t$ l'ensemble des courbes qui joignent les bords de l'anneau, et $\G_s$ la famille
des courbes de $A$ qui s\'eparent les composantes de bord. On d\'efinit alors
$$\mod_{\sup} (A,\SSS)= \frac{1}{\mod_2(\G_t,\SSS)}\quad et \quad \mod_{\inf} (A,\SSS)= \mod_2(\G_s,\SSS)\,.$$
De plus, si on note $E$ et $F$ les composantes du bord de $A$, on d\'efinit aussi $\hmod_{\sup}(A,\SSS)=1/\hmod_2(E,F,\SSS)$.
Lorsque $\SSS$ est planaire, on peut consid\'erer les courbes $\hG_s$ de $\NNN(\SSS)$ qui s\'eparent $\SSS(E)$
et $\SSS(F)$; on pose alors $\hmod_{\inf}(A,\SSS)=\hmod_2(\hG_s,\SSS)$.\end{defi}

\subsection{Comparaisons des modules}

Pour terminer ce paragraphe, nous donnons quelques propri\'et\'es des modules 
combinatoires qui sont d\^us essentiellement \`a J.W.\,Cannon, W.\,Floyd et
W.\,Parry (voir \cite{cfp1} par exemple).

\medskip

Si  deux recouvrements sont d'indice de recouvrement fini, alors les modules associ\'es sont comparables\,:

\REFLEM{bddolap} Si $\SSS$ et $\SSS'$ sont deux recouvrements de $X$ d'indice de recouvrement  fini $K$, alors, pour toute famille
de courbes $\G$ et tout $Q\ge 1$, on a 
$$\mod_Q(\G,\SSS)\asymp\mod_Q(\G,\SSS')$$ o\`u les constantes implicites ne d\'ependent que de
$K$ et de $Q$.
\ENDLEM

\Preuve 
Si $\rho$ est une m\'etrique admissible pour $\SSS$, on
d\'efinit $\si:\SSS'\to\R_+$ par $$\si(s')=\max\{\rho(s),\ s\in\SSS,\ s\cap s'\ne \emty\}\,.$$
Du coup, si $\g$ est une courbe de $\G$, alors
$$\ell_{\si}(\g,\SSS') \ge (1/K)\sum_{s'\in\SSS'(\g)}\sum_{s\in\SSS(s')}\rho(s)\ge (1/K)\sum_{s\in\SSS(\g)}\rho(s)$$
donc $\ell_{\si}(\g,\SSS') \ge(1/K)\ell_{\rho}(\g,\SSS)$. D'autre part,
$$V_Q(\si,\SSS')\le \sum_{s'\in\SSS'}\sum_{s\in\SSS(s')}\rho(s)^Q\le K\sum_{s\in\SSS}\rho(s)^Q$$
donc $V_Q(\si,\SSS')\le K V_Q(\rho,\SSS)$ et on en d\'eduit $\mod_Q(\G,\SSS')\lesssim\mod_Q(\G,\SSS)$.
On \'etablit le lemme en inversant les r\^oles de $\SSS$ et $\SSS'$.\endp

\begin{defi}[tuiles, toiture] Une {\sl tuile} est un compact connexe et une  {\sl toiture} est un recouvrement
localement fini par tuiles.\end{defi}

\REFPROP{tvb} Soit $X$ un espace m\'etrique propre  connexe, localement connexe par arcs 
recouvert par une toiture
$\SSS$ de valence born\'ee. Alors pour tout condensateur $(E,F)$, on a
$$
\hmod_Q(E,F,\SSS)\asymp\mod_Q(E,F,\SSS)$$ o\`u les constantes implicites ne d\'ependent que de
 la valence.
\ENDPROP

\begin{rema} Lorsque $X$ est une surface et $\SSS$ est planaire, alors la d\'emonstration ci-dessous
montre que $\hmod_{\inf}(A,\SSS)\asymp\mod_{\inf}(A,\SSS)$.\end{rema}
 
\Preuve
Tout d'abord, notons que sous les conditions de cette proposition, un ouvert connexe de $X$ est aussi connexe par arcs.

On note $\hG$ la famille de courbes dans $\NNN(\SSS)$ qui joignent $\SSS(E)$ et $\SSS(F)$. Si $\rho$ est optimale pour
$\mod_Q(E,F,\SSS)$, on d\'efinit $$\hrho(s)=\sum_{s'\in B_1(s,1)}\rho(s).$$

Soit $\hg\in\hG$. Notons d'une part $L$ la r\'eunion des pi\`eces dans $X$ qui n'intersectent pas $\SSS(\hg)$ et d'autre part
$K$ la r\'eunion des pi\`eces qui forment le support de $\hg$. Alors $L$ est ferm\'e car $\SSS$ est localement fini, 
$K$ est un continuum
car $\SSS$ est une toiture, et $K\subset X\setminus L$. Du coup, la composante connexe de $X\setminus L$ qui contient 
$K$ contient une courbe $\g\in\G(E,F)$. Cette courbe se trouve dans le $1$-voisinage de $K$ dans 
$\NNN(\SSS)$.
Donc  $\ell_{\hrho}(\hg)\ge \ell_\rho(\g)$.

D'autre part, l'in\'egalit\'e de H\"older implique $$\hrho(s)^Q\le (K+1)^{Q-1}\sum_{s'\in B(s,1)}\rho(s')^Q$$
donc $$V_Q(\hrho)\le (K+1)^{Q-1}\sum_s \sum_{s'\in B(s,1)}\rho(s')^Q\le (K+1)^Q V_Q(\rho)\,,$$
puisque chaque pi\`ece $s'$ appara\^{i}t dans au plus $(K+1)$ boules distinctes de rayon $1$.\endp

\begin{rema} Si $\SSS$ est un pavage par domaines de Jordan 
alors on obtient 
$\mod_{\sup}=\hmod_{\sup}$ et $\mod_{\inf}=\hmod_{\inf}$. En effet, \`a chaque courbe $\hg$
dans le graphe va correspondre une courbe de la surface $\g$ telle que $\SSS(\g)=\SSS(\hg)$.\end{rema}

Nous donnons une nouvelle d\'emonstration du lemme suivant ({\it cf.} \cite{cfp1}).

\begin{lemm}\label{infsup} Si $A$ est un anneau d'une surface $X$ recouverte par $\SSS$, alors 
$$\mod_{\inf}(A,\SSS)\le\mod_{\sup}(A,\SSS)\,.$$ \end{lemm}

\Preuve Il suffit de trouver une m\'etrique admissible $\rho$ pour laquelle 
$\mod_2(\G_s,\rho,\SSS)\le 1/\mod_2(\G_t,\rho,\SSS)$. 
Dans ce sens, nous consid\'erons la m\'etrique
extr\'emale qui nous fournit $\mod_{\sup}(A,\SSS)$ avec $L_{\rho}(\G_t,\SSS)=1$. 
D'apr\`es le Lemme \ref{lbeur}, il existe
$\G_0$ et des coefficients $(\la_\g)\in\R^{\G_0}$ tels que, pour toute $s\in\supp\,\rho$, 
$$\rho(s) = (1/2)\sum_{\g\in\G_0(s)}\la_\g\,.$$

On remarque qu'une courbe $\hat{\g}\in\G_s$ doit couper chaque courbe de $\G_0$ par le 
\th de Jordan, donc la $\rho$-longueur
de $\hat{\g}$ doit \^etre au moins $(1/2)\sum_{\G_0} \la_\g$ donc 
$$L_\rho(\G_s,\SSS)\ge (1/2)\sum_{\G_0} \la_\g\,.$$
D'autre part, on a 
$$\begin{array}{ll} \dis\sum_{\SSS}\rho(s)^2 & = (1/2)\dis\sum_{\SSS}\rho(s)\dis\sum_{\g\in\G_0(s)} \la_\g\\ & \\
& = (1/2)\dis\sum_{\g\in\G_0}\la_\g\dis\sum_{\SSS(\g)}\rho(s) \\ & \\
& = (1/2)\dis\sum_{\g\in\G_0}\la_\g\\ & \\
& \le L_\rho(\G_s)\,.\end{array}$$
Par cons\'equent, on obtient 
$$\frac{V_2(\rho,\SSS)}{L_{\rho}(\G_s,\SSS)^2}\le \frac{1}{V_2(\rho,\SSS)}$$ soit
$$\mod_{\inf}(A,\SSS)\le \mod_2(\G_s,\rho,\SSS)\le \frac{1}{\mod_2(\G_t,\rho,\SSS)}=\mod_{\sup}(A,\SSS)\,.$$
\endp

\begin{rema} Dans \cite{cfp2}, les auteurs vont jusqu'\`a montrer que dans le cas d'une toiture
de valence born\'ee, on a $\mod_{\inf}=\hmod_{\sup}$. Une version plus faible de ce r\'esultat sera \'etablie
au Corollaire \ref{inegtr}.\end{rema}

\section{Th\'eor\`eme de Riemann combinatoire et empilements de cercles}

Notre objectif ici est de m\'elanger les approches de J.W.\,Cannon et de M.\,Bonk et B.\,Kleiner pour
retrouver l'essentiel de leurs r\'esultats. On commence par une version faible du \th de Riemann combinatoire
qui nous servira ensuite \`a obtenir des param\'etrages quasisym\'etriques de sph\`eres topologiques soumises
\`a des conditions analytiques.

\subsection{Th\'eor\`eme de Riemann combinatoire revisit\'e}

On montre une version faible du \th de Riemann combinatoire (\Th \ref{comp}) en utilisant les empilements
de cercles. Notre approche en redonne une d\'emonstration succincte
dans un cadre simplifi\'e. 

\subsubsection{Le \th de Riemann combinatoire}
Rappelons l'enonc\'e du th\'eor\`eme original de J.W.\,Cannon. Pour cela on introduit quelques notions
suppl\'ementaires.

\begin{defi}[approximation s\'eparante] C'est une suite de recouvrements localement finis 
$(\SSS_n)_n$ dont la
maille tend vers $0$ telle que \ben
\item pour tout anneau $A$ , il existe $m>0$ et $n_0$ telle que, pour $n\ge n_0$,
on ait $$\mod_{\sup} (A,\SSS_n)\ge m\,;$$
\item pour tout $x\in X$, tout voisinage $V$ de $x$ et tout $m>0$, il existe un anneau $A\subset V$ qui s\'epare
$x$ de $X\setminus V$ et $n_0$ tels que $$\mod_{\sup}(A,\SSS_n)\ge m $$ pour $n\ge n_0$.\een 
\end{defi}

\begin{defi}[approximation conforme] C'est une  suite de recouvrements  $(\SSS_n)_n$ dont la
maille tend vers $0$ telle que \ben
\item la suite est une approximation s\'eparante\,;
\item il existe $K\ge 1$ telle que, pour tout anneau $A$ bord\'e par des courbes de Jordan, 
il existe $m>0$ et $n_0$ telles que, pour $n\ge n_0$, $$\mod_{\sup}(A,\SSS_n),\mod_{\inf}(A,\SSS_n)\in [m/K,Km].$$
\een \end{defi}

Avec ces notations, on a

\begin{theo}[de Riemann combinatoire \cite{ca}]\label{thm:rmt} Si $(\SSS_n)$ est une approximation conforme de toitures sur 
une surface topologique, alors celle-ci admet une structure complexe telle que les modules analytiques des anneaux
soient comparables \`a leurs modules combinatoires pour $n$ assez grand {\it i.e.}, il existe $K'$ telle que pour tout anneau $A$,
il existe $n_0=n_0(A)$ et $m>0$ telles que si $n\ge n_0$ alors 
$$\mod_{\sup}(A,\SSS_n),\mod_{\inf}(A,\SSS_n),\mod_2 A\in [m/K',K'm].$$\end{theo}

\subsubsection{Une version sph\'erique}
Soit $X$ une surface topologique hom\'eomorphe \`a $\SS^2$ que l'on suppose m\'etris\'ee. 
Soit $(\SSS_n)_n$ une suite
de toitures
dont la maille tend vers $0$ et telle que le nerf $\NNN(\SSS_n)$
de $\SSS_n$ contient une triangulation $\TTT_n$ de $\SS^2$ \'equivalente au sens de la d\'efinition \ref{geq} pour chaque $n$. 
Observons que l'existence d'un recouvrement dont le nerf est une triangulation provient du fait que $X$ est de dimension
topologique $2$, voir \cite{nag}. 

On remarque que puisque les triangulations sont planaires,  
on peut aussi consid\'erer $\hmod_{\inf}(\cdot,\TTT_n)$. 
Cependant, $\NNN(\SSS_n)$ peut ne pas \^etre planaire\,; c'est pourquoi on d\'efinit 
$\hmod_{\inf}(\cdot,\SSS_n)=\hmod_{\inf}(\cdot,\TTT_n)$ par commodit\'e.

On se fixe trois points
distincts $a_1,a_2,a_3\in X$.
Pour chaque $n\ge 0$ et chaque $s\in\SSS_n$, on se donne $p(s)\in s$ et on note $\PPP_n$ l'ensemble des ces points.
 On les choisit pour que $p(s)\ne p(s')$ si $s\ne s'$ et pour  que $\{a_j\}\subset \PPP_n$.
 
Pour chaque $n$, on consid\`ere l'empilement $\EEE_n=\EEE(\TTT_n)$ donn\'e par le \Th \ref{andreev}
tel que le centre du disque associ\'e \`a $p^{-1}(a_1)$ soit $0$, celui associ\'e \`a $p^{-1}(a_2)$ soit $1$ et celui
associ\'e \`a $p^{-1}(a_3)$ soit le point \`a l'infini. 

On d\'efinit l'application $\phi_n: \PPP_n\to \cbar$ qui \`a $p(s)$ associe le centre $c(s)$ du cercle
correspondant.

\begin{theo}\label{comp} Si la suite de toitures $(\SSS_n)$ est de valence born\'ee et 
s\'eparante 
alors on peut extraire une sous-suite $(n_k)$ telle 
que $(\phi_{n_k})$ converge uniform\'ement vers un \homeo $\phi:X\to\cbar$ 
tel que, pour tout anneau $A\subset X$,
$$\mod_{\sup}(A,\SSS_{n_k})\asymp\mod\,\phi (A)\asymp\mod_{\inf}(A,\SSS_{n_k})$$
pour $k$ assez grand. Autrement dit, la sous-suite est 
conforme.
\end{theo}

Dans la suite, il sera pratique d'\'etendre $\phi_n$ aux parties de $X$\,: \'etant donn\'e un ensemble $E\subset X$, 
on peut d\'efinir $\phi_n(E)$ comme la r\'eunion des tuiles du pavage associ\'e \`a $\EEE_n$ qui correspondent
\`a $\SSS_n(E)$. 
On d\'efinit de mani\`ere similaire $\phi_n^{-1}$. Notons que, $\SSS_n$ \'etant une toiture
et $\TTT_n\subset\NNN(\SSS_n)$, $\phi_n^{-1}$ transforme continua en continua.
Quant \`a $\phi_n$, on s'appuiera sur la Proposition \ref{mgr} (2).

La traduction de la s\'eparation 
se fera \`a l'aide du Corollaire \ref{inegtr} qui suit. 

Le fait de travailler avec de {\it v\'eritables} disques peut en fait \^etre all\'eg\'e. En g\'en\'eral, le
point important est de pouvoir contr\^oler leur masse par leur diam\`etre, ce qui nous motive \`a donner
la d\'efinition suivante.

\begin{defi}[Rondeur] Si $U\subset X$ est un ensemble born\'e et $x\in U$, alors on d\'efinit la {\sl rondeur
de $U$ par rapport \`a $x$} comme 
$$\rond(U,x) =\inf \{R/r,\ B(x,r)\subset U\subset B(x,R)\}\,.$$
Par $\rond(U)$, on entendra la plus petite valeur de $\rond(U,x)$ quand $x$ parcourt $U$.
De plus, si $K\ge 1$, on dira que $U$ est $K$-rond si $\rond (U)\le K$.\end{defi}

{\noindent\bf Pavage associ\'e \`a un empilement de cercles.}
Si $D_1$,  $D_2$,  $D_3$ sont trois disques ferm\'es d'int\'erieurs deux \`a deux disjoints de $\DD$ et sont tangents
deux \`a deux, on d\'esigne par $\oo$ une composante du compl\'ementaire 
de ces disques bord\'es par trois
arcs de cercles que l'on appellera {\sl interstice}. 
Cet ouvert $\oo$ est hyperbolique avec trois points marqu\'es, et il existe un unique
``centre'' {\it i.e.}, un point $c$ telle que la mesure harmonique de chaque ar\^ete vue de $c$ est exactement
$1/3$. 
D\'ecoupons $\oo$ en trois parties en consid\'erant les rayons g\'eod\'esiques issus de $c$ qui aboutissent
\`a chaque point marqu\'e. 
On associe alors aux disques la partie
qui lui est la plus proche. 
On d\'efinit alors $\hat{D}_j=s_j$ comme la r\'eunion de $D_j$ avec les portions  rajout\'ees par le processus 
pr\'ec\'edent \`a chaque point de tangence. Les rayons hyperboliques bordant $\hat{D}_j$ forment
ses {\sl c\^ot\'es}.

Si $\EEE$ est un empilement de cercles de $\cbar$ de diam\`etre plus petit que $\pi$, alors
on peut associer un pavage $\SSS=\SSS(\EEE)$ comme ci-dessus, et \'etudier les modules combinatoires associ\'es.
Le point crucial de cette construction est que $\rond (\hat{D}_j,c_j)\le K$ o\`u $c_j$ est le centre du disque
$D_j$ et $K$ est une constante universelle (par le \th de distorsion de Koebe).
Du coup $\SSS(\EEE)$ est un quasi-empilement de $\cbar$ au sens de la d\'efinition \ref{def:qemp}. 
On dira par la suite que $\SSS(\EEE)$ est la toiture
associ\'ee \`a l'empilement de cercles $\EEE$.

Dans le cadre des pavages associ\'es \`a des empilements de cercles, 
les diff\'erentes notions de modules combinatoires co\"{\i}ncident.

\begin{defi}[anneau empil\'e] Soit $\SSS$ la toiture associ\'ee \`a un empilement de cercles de $\cbar$.
Un anneau $A$ est dit empil\'e s'il est r\'eunion de tuiles de $\SSS$ et si chaque composante du compl\'ementaire
contient au moins un disque de l'empilement.\end{defi}

Dans le cas qui nous int\'eresse, le Lemme \ref{pongen} s'\'enonce ainsi.

\begin{lemm}[du pont]\label{linegtr} Soient $\SSS$ la toiture associ\'ee \`a un empilement de cercles $\EEE$ 
de $\cbar$,  et $\G$ une
famille de courbes. Si pour tout disque $D$ de $\EEE$ et toute courbe $\g\in\G$ qui traverse sa tuile
associ\'ee $\hat{D}$, 
on a $\ell(\g\cap 2KD)\ge \kappa \diam D$, o\`u $K$
est la rondeur de la toiture associ\'ee,
alors il existe une constante $C=C(\kappa)>1$ telle que
$$\mod_2\G\le C\cdot\mod_2(\G,\SSS).$$\end{lemm}

On en d\'eduit alors une version non-asymptotique de la Proposition \ref{disoka}\,:

\begin{coro}\label{inegtr} Si $\EEE$ est un empilement normalis\'e de valence born\'ee et 
si $A$ est un anneau empil\'e 
alors
$$\mod_{\sup} (A,\SSS)\asymp\mod\,A\asymp\mod_{\inf} (A,\SSS),$$
o\`u les constantes implicites ne d\'ependent que de la valence.\end{coro}

\Preuve On consid\`ere les familles de courbes $\G_t$ et $\G_s$ qui relient 
et s\'eparent les composantes de bord.
Pour appliquer le Lemme  \ref{linegtr}, il suffit de montrer l'existence d'une
constante $\kappa >0$ qui ne d\'epend que de la valence telle  que si $\g$
est l'une de ces courbes qui traverse une pi\`ece $s$, alors $\diam\g \ge\kappa \diam s$. 

Si $\g\in\G_s$, alors  $\g$ doit traverser la r\'eunion de ses voisines
et contourner au moins un disque\,;
or  le lemme du collier montre que ces  pi\`eces ont toute une taille comparable,
ce qui montre l'existence de $\kappa>0$ dans ce cas.

Si $\g\in\G_t$, alors $\g$ doit \^etre de diam\`etre au moins celui d'un
c\^ot\'e de $s$ ou de ses voisines. Comme les rayons des cercles bordant 
un interstice sont comparables, un argument de compacit\'e
montre que le diam\`etre des c\^ot\'es est aussi comparable \`a ces rayons.

Donc les  hypoth\`eses du Lemme \ref{linegtr} sont satisfaites pour
un $\kappa>0$ qui ne d\'epend que de la valence. 
De plus, le Lemme \ref{infsup} 
affirme que $\mod_{\inf}\le \mod_{\sup}$, 
donc on obtient 
$$\mod_{\sup} (A,\SSS)\lesssim\mod\,A\lesssim\mod_{\inf} (A,\SSS)\le\mod_{\sup} (A,\SSS).$$
\endp

{\noindent\sc D\'emonstration du \Th \ref{comp}.} On suppose que $(\SSS_n)$ est 
s\'eparante et de valence born\'ee. 
Soit $A\subset X$ un anneau dont les composantes du compl\'ementaire sont $E$ et $F$.
Pour $n$ assez grand, on associe le condensateur $(\hat{E}_n,\hat{F}_n)$ 
dans $\TTT_n$ par la Proposition  \ref{mgr} (2),
et l'anneau empil\'e $A_n$ qui lui correspond sur $\cbar$.
On a donc par les Propositions \ref{tvb} et \ref{mgr}
donc 
\begin{equation}\label{eqn:trans}
\mod_{\sup}(A,\SSS_n)\asymp\hmod_{\sup}(A,\SSS_n)\asymp \hmod_{\sup}(A_n,\SSS(\EEE(\TTT_n)))= 
\mod_{\sup}(A_n,\SSS(\EEE(\TTT_n)))\end{equation}

On obtient aussi $\mod_{\inf}(A_n,\SSS(\EEE(\TTT_n))) \lesssim \mod_{\inf}(A,\SSS_n)$,
car on peut associer \`a une courbe $\g$ de $\G_s(A_n)$ une courbe de $\G_s(A)$
dont le support est dans le $1$-voisinage de $\g$ dans $\NNN(\SSS_n)$, 
{\it cf.} la d\'emonstration de la Proposition \ref{tvb}. 

Montrons d'abord que la suite $(\phi_n)_n$ est \'equicontinue en chaque point.
On se fixe $\ep>0$ et $m>0$, et on consid\`ere $x\in X$. On peut supposer que $x$ est diff\'erent de
$\{a_1,a_2\}$. Soit alors $V$ un voisinage de $x$ disjoint de $\{a_1,a_2\}$. Comme l'approximation est
s\'eparante, il existe un anneau $A\subset V$ (bord\'e par des courbes de Jordan) qui
s\'epare $x$ de $X\setminus V$ tel que, pour $n$ assez grand, on ait 
$\hmod_{\sup} (A_n,\TTT_n)\ge m$, o\`u $A_n$ est l'anneau empil\'e associ\'e. 
Le Corollaire \ref{inegtr} implique $$\mod\,A_n \gtrsim \hmod_{\sup} (A_n,\TTT_n)\ge m\,.$$
Or pour $n$ assez grand, une des \ccs de $\cbar\setminus A_n$ contient $\{0,1\}$, donc le diam\`etre sph\'erique
de l'autre composante est born\'e par une constante $C(m)$ qui tend vers $0$ quand $m\to\infty$. 
On se fixe $m$ assez grand
pour que $C(m)\le\ep$.

Comme la maille des $(\SSS_n)$ tend vers $0$, il existe un voisinage $W\subset V$ de $x$ disjoint de tous les
$A_n$ et s\'epar\'e de $\{a_1,a_2\}$. Du coup, pour tout $n$ assez grand, pour tous $s,s'\subset W$, on a
$d(\phi_n(p(s)),\phi_n(p(s')))\le C(m)\le\ep$.
Ceci implique que la suite $(\phi_n)$ est \'equicontinue, donc on peut extraire une sous-suite $(\phi_{n_k})_k$
convergente de limite une application continue $\phi:X\to\cbar$. Nous renvoyons \`a l'Appendice
\ref{apphaus} pour plus de d\'etails sur cette convergence. 

On proc\`ede en deux \'etapes pour montrer que $\phi$ est injective. Tout d'abord, si $x\not\in\{a_j\}$,
alors on peut trouver un anneau $A$ bord\'e par des courbes de Jordan qui s\'epare $\{x,a_1\}$ de $\{a_2,a_3\}$ 
de module combinatoire uniform\'ement minor\'e. Par le Corollaire \ref{inegtr}, on en d\'eduit que les modules
analytiques des anneaux empil\'es associ\'es 
sont aussi uniform\'ement minor\'es, donc
$\phi(x)\notin\{1,\infty\}$, sinon le module devrait tendre vers z\'ero. 
Par permutation cyclique, on trouve que $\phi(x)\notin\{0,1,\infty\}$.

Dans un second temps, on consid\`ere $x\ne y$. On peut supposer que $x,y,a_1,a_2$ sont quatre points distincts.
En consid\'erant un anneau qui s\'epare $\{x,a_1\}$ de $\{y,a_2\}$, on montre que $\phi(x)\ne\phi(y)$ par le m\^eme
argument.

Par cons\'equent, $\phi$ est un \homeo de $X$ qui est une sph\`ere topologique  sur son image.
On en d\'eduit que $\phi:X\to\cbar$ est un hom\'eomorphisme.

\bigskip

La comparaison des diff\'erents modules associ\'es \`a un anneau suit aussi du Corollaire \ref{inegtr} (voir aussi le \Th 7.1 dans \cite{ca}).
En effet, si $A\subset \cbar$ et $(\EEE_{n_k})$ est la suite d'empilements associ\'ee \`a $(\TTT_{n_k})$, 
alors $(\SSS(\EEE_{n_k}))$ est conforme. 
En effet, si $A\subset\cbar$ est un anneau bord\'e par des courbes de Jordan, on note $A_k$ l'anneau
empil\'e associ\'e \`a $\EEE_{n_k}$. 
Comme la maille des recouvrements $\SSS_n$ tend vers $0$ et puisque la suite $(\phi_{n_k})$
est convergente, on en d\'eduit que la maille des $(\SSS(\EEE_{n_k}))$ tend aussi vers $0$. Du coup,
 $\mod\,A_k\to \mod\,A$, et 
les hypoth\`eses du Corollaire
\ref{inegtr} sont v\'erifi\'ees pour $\G_s(A_k)$ et $\G_t(A_k)$. Donc, pour $k$ assez grand,
on obtient ainsi 
$$\mod_{\sup} (A,\SSS(\EEE_{n_k})) \asymp\mod_{\inf} (A,\SSS(\EEE_{n_k}))\asymp \mod\,A_k\asymp \mod\,A\,.$$

Consid\'erons maintenant un anneau $A\subset X$ bord\'e par des courbes de Jordan. 
L'anneau $\phi(A)$ s'obtient donc comme limite d'anneaux empil\'es $A_k$ 
donc $\mod\, A_k$ tend vers $\mod\, \phi(A)$. Du coup, puisque la maille des empilements tend vers $0$,
$$\hmod_{\sup}(A,\SSS_{n_k})\asymp \hmod_{\sup}(\phi_{n_k}(A),\SSS(\EEE_{n_k}))\asymp\mod\, \phi(A)$$
pour $k$ assez grand. La Proposition \ref{tvb} montre que $\mod_{\sup}(A,\SSS_{n_k})\asymp\mod\, \phi(A)$, 
{\it cf.} (\ref{eqn:trans}).
Par ailleurs, $$\mod_{\sup}(A,\SSS_{n_k})\ge \mod_{\inf}(A,\SSS_{n_k})\gtrsim 
\mod_{\inf}(A_{n_k},\TTT_{n_k})\asymp\mod\, \phi(A)\,.$$
\endp

\begin{rema} L'\'equicontinuit\'e de $(\phi_n)_n$ \'etablie dans la d\'emonstration montre notamment que, 
pour tout anneau $A$, il existe
une constante $K=K(A)$ telle que les modules combinatoires de $A$ soient $K$-comparables
entre eux. On peut  se demander si cela n'implique pas que la suite soit automatiquement
conforme\,: la r\'eponse est n\'egative. En effet, soit $(\EEE_n)_n$ une suite d'empilements
de valence born\'ee dont la maille tend vers $0$ et de pavage associ\'e $(\PPP_n)$\,; 
consid\'erons aussi un \homeo $h:\cbar\to\cbar$ non quasiconforme, et posons 
$\SSS_{2n}=\PPP_n$ et $\SSS_{2n+1}=h^{-1}(\PPP_n)$. La suite $(\SSS_n)$ est bien
s\'eparante, mais elle n'est pas conforme puisque $h$ n'est pas quasiconforme.

Ceci montre que sous nos hypoth\`eses, deux structures complexes construites par le \th ne
sont pas forc\'ement dans la m\^eme classe quasiconforme (si $\Phi$ et $\Psi$ sont deux
\homeos limites, $\Phi\circ\Psi^{-1}$ n'est pas forc\'ement quasiconforme). Une modification
de l'exemple ci-dessus permet aussi de montrer que le fait que toute limite est dans la m\^eme
classe quasiconforme n'implique pas que la suite est conforme\,: consid\'erons maintenant
une suite d'\homeos \qcs $(h_k)$ qui tend uniform\'ement vers l'identit\'e avec une distorsion
qui diverge. On construit une suite de pavages en utilisant $h_k^{-1}(\PPP_n)$, $k,n\in\N$\,:
cette suite sera bien 
s\'eparante, les \homeos limites seront les $h_k$ ou 
l'identit\'e, mais la suite n'est pas conforme car la suite $(h_k)$ n'est pas uniform\'ement 
quasiconforme.\end{rema}

\medskip

On en d\'eduit une version faible du th\'eor\`eme de Riemann combinatoire. 

\begin{coro}\label{oricomp} Si la suite $(\SSS_n)$ est une suite de toitures de valence born\'ee et conforme, chaque $\SSS_n$ contenant une triangulation \'equivalente $\TTT_n$, alors 
 on peut extraire une sous-suite $(n_k)$ telle 
que $(\phi_{n_k})$
converge uniform\'ement vers un \homeo $\phi:X\to\cbar$ tel que 
$$\mod_{\sup}(A,\SSS_n)\asymp\mod\,\phi (A)\asymp\mod_{\inf}(A,\SSS_n)$$
pour $n$ assez grand. De plus, si $\phi$ et $\psi$ sont deux \homeos limites, alors
$\phi\circ\psi^{-1}$ est $K$-\qc o\`u $K$ ne d\'epend que des donn\'ees.

R\'eciproquement, si $(\phi_n)$ est convergente vers un \homeo alors $(\SSS_n)$ est conforme.
\end{coro}

\medskip
\Preuve 
Le \Th \ref{comp} s'applique pour montrer la convergence d'une sous-suite $(\phi_{n_k})_k$ vers
un \homeo $\phi$. Conjugu\'e
\`a la conformit\'e de la suite,  ce \th fournit aussi les estimations entre les modules,
et permet de montrer que si $\phi$ et $\psi$ sont deux limites, alors  $\phi\circ\psi^{-1}$ pr\'eserve les modules \`a une constante multiplicative uniforme pr\`es, autrement dit $\phi\circ\psi^{-1}$ est $K$-\qc   o\`u $K$ ne d\'epend que des donn\'ees ({\it cf.} appendice \ref{tgf}).

\medskip

Quant \`a la r\'eciproque, elle d\'ecoule du m\^eme argument que dans la d\'emonstration du \Th \ref{comp}.\endp

\subsection{Param\'etrage quasisym\'etrique de sph\`eres topologiques}

Dans ce paragraphe, on exhibe des conditions (analytiques) sur une sph\`ere topologique qui permettent
de conclure qu'elle est quasisym\'etrique \`a la sph\`ere de Riemann. Notre point de d\'epart sera
toujours la donn\'ee d'une suite de quasi-empilements par tuiles
de valence born\'ee. Nous verrons comment,  avec l'approche du paragraphe pr\'ec\'edent,
\'etablir des r\'esultats de M.\,Bonk et B.\,Kleiner.

On se r\'ef\`ere  de mani\`ere essentielle \`a l'appendice \ref{tgf} pour les d\'efinitions, les notations et les
propri\'et\'es des objets analytiques utilis\'ees dans ce paragraphe. 

Le premier r\'esultat nous servira pour conclure que les limites des $(\phi_n)$ construites au paragraphe pr\'ec\'edent
sont quasisym\'etriques.

\begin{prop}\label{regqs}
Soient $X$ un espace doublant 
et $(\SSS_n)$ une suite de quasi-empilements qui v\'erifie les hypoth\`eses 
du \Th \ref{comp}. 
Les propositions suivantes sont \'equivalentes.
\ben \item Il existe $\eta$ telle que tout \homeo limite $\phi$ est $\eta$-quasisym\'etrique.
\item Il existe un \homeo 
$\psi_1:\R_+\setminus\{0\}\to\R_+\setminus\{0\}$
d\'ecroissant
tel que 
pour tout condensateur $(E,F)$, on ait $\mod_2(E,F,\SSS_n)\ge \psi_1(\De(E,F))$ pour tout
$n$ assez grand.
\item L'espace $X$ est  lin\'eairement localement connexe et il existe un \homeo 
$\psi_2:\R_+\setminus\{0\}\to\R_+\setminus\{0\}$
d\'ecroissant 
tel que,  pour tout condensateur $(E,F)$, on a $\mod_2(E,F,\SSS_n)\le \psi_2(\De(E,F))$ pour tout
$n$ assez grand.
\een
De plus, 
la suite  $(\SSS_n)_n$ est conforme.
\end{prop}

\Preuve 
Une application quasisym\'etrique contr\^ole le rapport de distances, en particulier les distances relatives $\De(E,F)$
({\it cf.} Lemme 3.2 dans \cite{bk1}). Donc si $\phi$ est $\eta$-quasisym\'etrique, alors il existe des \homeos 
$\chi_1,\chi_2:\R_+\to\R_+$ ne d\'ependant que de $\eta$ tels que $$\chi_1(\De (E,F))\ge \De(\phi E,\phi F)\ge \chi_2(\De (E,F))\,.$$

Sur $\cbar$, il existe des \homeos $\hat{\psi}_1,\hat{\psi}_2:\R_+\setminus\{0\}\to\R_+\setminus\{0\}$
d\'ecroisssantes vers $0$ telles que $$  \hat{\psi}_2(\De (E, F))\ge \mod_2(E,F)\ge  \hat{\psi}_1(\De (E, F))\,.$$

On se fixe $(E,F)$ et on d\'efinit
$$m_s=\limsup \mod_2(E,F,\SSS_n)\quad\hbox{et}\quad m_i=\liminf \mod_2(E,F,\SSS_n)\,.$$

Montrons que 
\begin{equation}\label{eqn:mis} m_s\lesssim m_i
\end{equation} 
o\`u les constantes implicites ne d\'ependent pas
du condensateur. On se donne les sous-suites $(\SSS_{n_k,s})$ et $(\SSS_{n_k,i})$
qui font tendre les modules vers $m_s$ et $m_i$. Le \Th \ref{comp} produit deux \homeos $\eta$-quasisym\'etriques
$\phi_s,\phi_i:X\to\cbar$ tels que $$m_s\asymp \mod_2(E,F,\SSS_{n_k,s})\asymp \mod_2( \phi_s(E),\phi_s(F))$$
et $$m_i\asymp \mod_2(E,F,\SSS_{n_k,i})\asymp \mod_2( \phi_i(E),\phi_i(F))\,.$$
Or $\phi_s\circ\phi_i^{-1}$ est un \homeo $\heta$-quasisym\'etrique de la sph\`ere 
de Riemann, o\`u $\heta(t)=\eta(1/\eta^{-1}(1/t))$.
Donc $\phi_s\circ\phi_i^{-1}$ quasipr\'eserve les $2$-modules. Par cons\'equent, $m_s\lesssim m_i$ \'etablissant (\ref{eqn:mis}).

Cela signifie qu'il suffit de  travailler avec une sous-suite $(\SSS_{n_k})$ pour laquelle
on a convergence des injections vers un \homeo $\phi$.

D'apr\`es le \Th \ref{comp}, 
$$\mod_2(E,F,\SSS_{n_k})\asymp\mod_2(\phi E,\phi F)$$ 
 pour $k$ assez grand,
donc, pour $n$ assez grand et en vertu de (\ref{eqn:mis})
$$\left\{\begin{array}{ll} \mod_2(E,F,\SSS_{n})\asymp \mod_2(\phi E,\phi F) &
\ge  \hat{\psi}_1(\De (\phi E,\phi  F))\ge \hat{\psi}_1\chi_1(\De (E, F))\\
 \mod_2(E,F,\SSS_{n})\asymp \mod_2(\phi E,\phi F) &
\le  \hat{\psi}_2(\De (\phi E,\phi  F))\le \hat{\psi}_2\chi_2(\De (E, F)),\end{array}\right.
$$
ce qui montre l'existence des applications $\psi_1$ et $\psi_2$ sous la condition que $\phi$ est quasisym\'etrique.
Notons que la sph\`ere est lin\'eairement localement connexe, et cette propri\'et\'e est pr\'eserv\'ee par applications
quasisym\'etriques (cf. \cite[Chap.\,15]{he}).

Pour la r\'eciproque, on note que $X$ \'etant doublant, il suffit de montrer que toute limite $\phi$ est uniform\'ement
faiblement quasisym\'etrique
(voir la fin de l'appendice \S\,\ref{tgf}, ou \cite[Thm 10.19]{he}).
On se donne $x\in X$ et $r>0$ et on note $L=L_{\phi}(x,r)$, $\ell=\ell_{\phi}(x,r)$, 
$A=D(\phi (x), L)\setminus D(\phi (x), \ell)$, o\`u $D(z,r)\subset\cbar$ est le disque
centr\'e en $z$ de rayon $r$ pour la m\'etrique sph\'erique, et 
$\G_t$ la famille des courbes qui joignent les deux composantes
de bord de $A$. Si $L/\ell\le 2$, alors, on n'a rien \`a montrer. Sinon, on remarque que 
$\De(A)\asymp(L-\ell)/\ell\asymp L/\ell$, o\`u $\De(A)$ d\'esigne la distance relative des deux
composantes de bord de $A$.

D'une part, on a $$0<c \le\psi_1(\De(\phi^{-1}(A)))\le \mod_2 (\G_t,\SSS_n)$$ par hypoth\`eses, 
et d'autre part, 
si $n$ est assez grand, 
$$\mod_2 (\G_t,\SSS_n)\lesssim 1/\log L/\ell\,.$$
En effet, il suffit de consid\'erer les poids $\rho_n(s)= \diam\phi_n(s)/\dist(\phi(x),\phi_n(s))$ pour
$s\in\SSS_n(A)$ et $\rho_n(s)=0$ sinon ({\it cf.} Lemme \ref{inegAR}, et la d\'emonstration du Corollaire \ref{bk10.4}).
Donc ces deux in\'egalit\'es montrent que le rapport $L/\ell$ est born\'e ind\'ependamment de $\phi$, $x$ et $r>0$, donc $\phi$ est 
uniform\'ement faiblement quasisym\'etrique.

Quant \`a l'autre condition, (3), on montre que $\phi^{-1}$ est quasisym\'etrique en inversant les r\^oles. 
Soit $C$ la constante de connexit\'e locale lin\'eaire de $X$. Soient $z,z_1,z_2\in\cbar$ tels que $|z-z_1|\le |z-z_2|$.
Si $|\phi^{-1}(z)-\phi^{-1}(z_1)|\le 2C^2 |\phi^{-1}(z)-\phi^{-1}(z_2)|$, nous n'avons rien \`a montrer.
Sinon, on peut supposer que $z,z_1,z_2$ sont loin de l'infini et, par continuit\'e uniforme que leurs images
par $\phi^{-1}$ sont loin de $a_3$. On  construit deux continua disjoints $E$ et $F$ sur $X$ tels que 
$$\phi^{-1}(z),\phi^{-1}(z_2)\in E\subset B(\phi^{-1}(z), C |\phi^{-1}(z)-\phi^{-1}(z_2)|)$$
et 
$$\phi^{-1}(z_1), a_3\in F\subset X\setminus B(\phi^{-1}(z), (1/C) |\phi^{-1}(z)-\phi^{-1}(z_1)|)\,.$$
On a donc $$\De(E,F)\ge (1/C^3)\frac{|\phi^{-1}(z)-\phi^{-1}(z_1)|}{ |\phi^{-1}(z)-\phi^{-1}(z_2)|}$$
et $\De(\phi(E),\phi(F))\le 1$. Par suite,
$$0<\hat{\psi}_1(1)\lesssim \hmod_2(E,F,\SSS_n)\le 
\psi_2\left((1/C^3)\frac{|\phi^{-1}(z)-\phi^{-1}(z_1)|}{ |\phi^{-1}(z)-\phi^{-1}(z_2)|}\right)\,.$$
Comme $\psi_2$ est d\'ecroissante, on en d\'eduit que $|\phi^{-1}(z)-\phi^{-1}(z_1)|\lesssim |\phi^{-1}(z)-\phi^{-1}(z_2)|$,
et donc que $\phi^{-1}$ est quasisym\'etrique.

\medskip

L'estimation (\ref{eqn:mis}) s'interpr\`ete par la conformit\'e de $(\SSS_n)$ gr\^ace \`a la conclusion
du \Th \ref{comp} sur les estimations des modules.
\endp

\begin{rema}\label{rmq:tuilevsboules} Dans \cite{bk1}, M.\,Bonk et B.\,Kleiner travaillent
 avec ce qu'ils appellent des
{\em $K$-approximations} $(\AAA_n)$ qui ne sont pas des toitures (voir \S\,4). 
Cependant, ce sont des quasi-empilements de valence born\'ee qui ont des propri\'et\'es suppl\'ementaires 
\`a nos suites de recouvrements\,; notamment, deux points d'une m\^eme pi\`ece $U\in\AAA_n$ peuvent \^etre joints
par une courbe dont le support est contenu dans un $K$-voisinage de $U$ dans le nerf de $\AAA_n$ et de diam\`etre
comparable \`a celui de $U$.
Donc, en agrandissant les pi\`eces de $\AAA_n$, on obtient  un quasi-empilement par tuiles $\SSS_n$
(a) de valence born\'ee dont le  nerf est \'equivalent \`a celui de $\AAA_n$ au sens de la d\'efinition \ref{geq}
et (b) d'indice de recouvrement fini. Par cons\'equent, les Propositions \ref{mgr} et \ref{tvb}, et le Lemme \ref{bddolap} 
impliquent
que les modules combinatoires $\hmod_*$ et $\mod_*$ 
sont tous comparables, et nous restons bien dans le cadre qu'ils traitent.

Notons que la construction de ces recouvrements est d\'elicate
et  repose sur le fait que $X$ est doublant et lin\'eairement
localement connexe.  M.\,Bonk et B.\,Kleiner d\'efinissent
leurs approximations en construisant une triangulation
plong\'ee dans $X$, voir \cite[\S\,6]{bk1}. 
\end{rema}

Nous pouvons directement en d\'eduire ce que M.\,Bonk et B.\,Kleiner appellent le r\'esultat principal
de leur article \cite{bk1} (\Th 11.1)\,:

\begin{coro}[M.\,Bonk \& B.\,Kleiner]\label{th11.1}  Soit $X$ un espace m\'etrique doublant et lin\'eairement localement connexe 
hom\'eomorphe \`a $\SS^2$.
On suppose que
l'on a une suite $(\SSS_n)$ de quasi-empilements par tuiles 
dont la maille tend vers z\'ero et 
dont les nerfs contiennent 
des triangulations \'equivalentes de la sph\`ere
de valence born\'ee. 

S'il existe des constantes $\la>1$ et $C>0$ telles que
$$\hmod_2(B,X\setminus\la B,\SSS_n)\le C$$ pour toute boule $B\subset X$ et pour $n$ assez grand, alors
$X$  est quasisym\'etrique \`a $\cbar$ et $(\SSS_n)$ est conforme.
\end{coro}

\Preuve Soient $x\in X$ et $r>0$ et consid\'erons, comme M.\,Bonk et B.\,Kleiner, 
les anneaux $B(x,\la^{-2j}r)\setminus B(x,\la^{-(2j+1)}r)$.
Par l'in\'egalit\'e de Gr\"otzsch \cite[Thm 4.2]{ah1} et nos hypoth\`eses, on obtient pour $n$ assez grand
$$\hmod_{\sup} (B(x,r)\setminus B(x,\la^{-(2k+1)}r),\SSS_n)\ge \sum_{j=0}^k \hmod_{\sup} (B(x,\la^{-2j}r)\setminus B(x,\la^{-(2j+1)}r),\SSS_n) \ge k/C\,.$$
Donc $(\SSS_n)$ v\'erifie la condition (2) d'approximation 
s\'eparante qui permet de montrer que la suite
$(\phi_n)$ est \'equicontinue. 
S'il y avait une fibre non triviale $F$, alors
on pourrait trouver une boule $B\subset X$ telle que les diam\`etres de $\phi(B)$ et de 
$\phi(X\setminus \la B)$ sont minor\'es, avec les intersections $B\cap F$ et $(X\setminus\la B)\cap F$ non vides, donc 
telle que
$\phi(B)\cap\phi(X\setminus \la B)\ne\emty$. 
Par le Corollaire \ref{inegtr}, ceci contredit $\hmod_2(B,X\setminus\la B,\SSS_n)\le C$ 
pour $n$ assez grand. Donc les conclusions du \Th \ref{comp} sont v\'erifi\'ees.

 La Proposition \ref{regqs} (3) s'applique, ce qui montre  que
la limite $\phi$ est quasisym\'etrique et $(\SSS_n)$ est conforme.
\endp

On en d\'eduit ({\it cf.} \Th 10.4 de \cite{bk1})\,:

\begin{coro}\label{bk10.4} 
Si $(X,\mu)$ est un espace m\'etrique $2$-Ahlfors r\'egulier lin\'eairement localement connexe hom\'eomorphe \`a $\SS^2$, et 
si $(\SSS_n)$ est une suite de 
quasi-empilements par tuiles
de valence born\'ee dont le nerf contient une triangulation \'equivalente et la maille tend vers $0$ alors 
$(\SSS_n)$ est conforme et 
$X$ est quasisym\'etrique \`a $\cbar$.\end{coro}

\Preuve Il suffit de v\'erifier que les hypoth\`eses du Corollaire \ref{th11.1} sont v\'erifi\'ees.
On se fixe $x\in X$, et $r<R=2^kr<\diam X$, o\`u $k\ge 1$. 
On d\'efinit la m\'etrique $\rho_n:\SSS_n(B(x,R)\setminus B(x,r))\to\R_+$
par $\rho_n(s)= \diam s/\mbox{dist\,} (x,s)$ ({\it cf.} Lemme \ref{inegAR}). Alors, en d\'ecoupant en couronnes diadiques
centr\'ees en $x$, pour toute $\g\in\hG_t$, on a 
$$\ell_{\rho_n}(\g) \ge (1/r)\dis \sum_{j=0}^{k-1} 2^{-j} \dis \sum \diam s
\gtrsim  \dis\sum_{j=1}^k 2^{-j} 2^j
 \gtrsim k  \gtrsim \log R/r.$$

De m\^eme, en utilisant la r\'egularit\'e de $X$, les propri\'et\'es de rondeur et de valence born\'ee, on a
$$\begin{array}{ll}
 \dis\sum \rho_n(s)^2
&  \lesssim  \dis\sum_{j=0}^{k-1} \dis\sum\dis\frac{(\diam s)^{2}}{(2^{j} r)^{2}}\\ &\\
&  \lesssim  \dis\sum_{j=0}^{k-1} \dis\frac{\mu(B(x,2^{j+1} r))}{(2^{j} r)^{2}}\\ &\\
&  \lesssim  \dis\sum_{j=0}^{k-1} \dis\frac{(2^{j+1}r)^{2}}{(2^{j} r)^{2}}\\ &\\
& \lesssim \log R/r\,.
\end{array}$$
Donc $\hmod_2(B(x,r),X\setminus B(x,R),\SSS_n)\lesssim 1/\log R/r$, et le Corollaire \ref{th11.1} s'applique.
\endp

\medskip

On peut aussi retrouver une version asymptotique du \Th 10.1 de M.\,Bonk et B.\,Kleiner dans \cite{bk1}\,:

\begin{theo}\label{bk10.1} Soit $X$ un espace m\'etrique $Q$-Loewner et $Q$-Ahlfors r\'egulier
hom\'eomorphe \`a $\SS^2$. 
On suppose que
l'on a une suite de 
quasi-empilements par tuiles de valence born\'ee dont la maille tend vers z\'ero et dont le nerf contient 
une triangulation \'equivalente de la sph\`ere. 
Alors 
$(\SSS_n)$ est conforme, 
$X$  est quasisym\'etrique \`a $\cbar$ et $Q=2$.
\end{theo} 

\Preuve 
La condition loewnesque permet de minorer les $Q$-modules de courbes\,:
elle nous fournit l'existence d'une application d\'ecroissante $\psi_1$ telle que, pour tout condensateur $(E,F)$
de $X$,
$$\mod_Q(E,F)\ge \psi_1(\De (E,F)).$$ 
Or, $Q\ge 2$ implique que 
$\mod_Q(\G,\SSS_n)\le\mod_2(\G,\SSS_n)$ en normalisant les m\'etriques par $L_\rho(\G)=1$ pour toute famille de courbes\,; donc
 on obtient avec le Lemme \ref{pongen}
$$\psi_1(\De (E,F))\le \mod_Q(E,F)\lesssim \mod_Q(E,F,\SSS_n) \le\mod_2(E,F,\SSS_n)\,. $$
De plus, le Corollaire \ref{inegtr} et le Lemme \ref{inegAR} impliquent l'existence d'une autre
 application d\'ecroissante $\psi_2$ telle que
$\mod_2(E',F',\SSS(\EEE_n))\lesssim \mod_2(E',F')\le \psi_2(\De (E',F'))$ pour tout condensateur $(E',F')$ de $\cbar$ et $n$ assez grand. 
En particulier, pour tout condensateur $(E,F)$ de $X$,
et tout $n$ assez grand, on a 
\begin{equation}\label{eqn:1} \psi_1(\De(E,F))\lesssim \psi_2(\De(\phi_n E,\phi_n F))\,.
\end{equation}

Nous allons montrer que la suite $(\phi_n^{-1})_n$ est \'equicontinue, et que toute limite produit
un \homeo de $\cbar$ sur $X$. Pour cela, on montre d'abord que la maille des empilements de cercles
tend vers $0$. 

Si ce n'est pas le cas, il existe $r_0>0$, une sous-suite $(n_k)$ et un \'el\'ement
$s_k\in\SSS_{n_k}$ tels que son disque associ\'e $D_k$ est de rayon au moins $r_0$. Quitte 
\`a extraire une sous-suite, on peut supposer que $(s_k)$ tend vers un point $x\in X$ diff\'erent de
$a_1,a_2$. 
Rappelons que chaque nerf $\NNN(\SSS_n)$ contient la triangulation \'equivalente $\TTT_n$\,;
il existe donc une constante $K$ telle que si $s,s'\in\SSS_n$ sont \`a distance combinatoire au moins
$K$ dans $\TTT_n$, alors leur distance est au moins $2$ dans $\NNN(\SSS_n)$. Soit $\g_k\subset\cbar$
un arc qui joint $0$ \`a $D_k$\,; 
notons $E_k\subset\cbar$ le continuum constitu\'e de $\SSS_{n_k}(\g_k)$ et $V_k$ la r\'eunion
des tuiles qui forment le $K$-voisinage de $E_k$ dans la triangulation $\TTT_{n_k}$. Pour $k$ assez
grand, $V_k$ est bord\'e par une courbe de Jordan puisque la triangulation est planaire, 
et quitte \`a modifier $\g_k$, on peut aussi
supposer que $1$ n'y appartient pas. On construit alors un autre continuum $F_k\subset \cbar\setminus
V_k$ form\'e de tuiles de $\SSS(\EEE_{n_k})$ qui relie $1$ \`a une pi\`ece $\hat{D_k'}=\phi_{n_k}(s_k')$, o\`u
$s_k'$ est \`a distance $2K$ de $s_k$ dans $\TTT_{n_k}$, de sorte que $\De(E_k,F_k)\ge \de_0 >0$, o\`u $\de_0$ est une constante
ind\'ependante de $k$\,; l'existence de $\de_0$ provient du fait que les disques $D_k,D_k'$ associ\'es
\`a $s_k$ et $s_k'$ sont \`a distance
combinatoire uniform\'ement born\'ee et de rayon minor\'e, {\it cf.} le lemme du collier.
Or $\De(\phi_{n_k}^{-1}(E_k),\phi_{n_k}^{-1}(F_k))$ tend vers $0$, car 
les diam\`etres de $\phi_{n_k}^{-1}(E_k)$ et $\phi_{n_k}^{-1}(F_k)$ 
sont minor\'es par $\min\{d_X(a_1,x),d_X(a_2,x)\}/2$,  la maille de $\SSS_{n_k}$ tend vers $0$
et $s_k$ et $s_k'$ sont \`a distance combinatoire au plus $2K$. Ceci contredit (\ref{eqn:1}).

Il s'ensuit que la suite d'empilements de cercles $(\EEE_n)_n$ a sa maille qui tend vers $0$,
donc le Corollaire \ref{inegtr} implique que la suite des toitures associ\'ees est conforme. La majoration (\ref{eqn:1}) implique facilement
l'\'equicontinuit\'e de $(\phi_n^{-1})$. Notons que les arguments de la d\'emonstration du 
\Th \ref{comp} permettent de montrer l'injectivit\'e de la limite, ainsi que la conformit\'e
de $(\SSS_n)$. La
Proposition \ref{regqs} implique alors que $\phi^{-1}$ est quasisym\'etrique.

Notons que $(\SSS_n)$ \'etant conforme, on a forc\'ement $Q=2$. Pour $Q>2$, on aurait $\lim_n \mod_Q(\cdot,\SSS_n)=0$.
En effet, si $(E,F)$ est un condensateur de $X$ et $\phi:X\to\cbar$ est l'\homeo donn\'e par le \Th \ref{comp}, 
alors on consid\`ere $\rho_n(s)=\diam\phi(s)$ pour $s\in\SSS_n$. Il vient, pour $\g\in\G(E,F)$,
$$\ell_{\rho_n}(\g,\SSS_n)=\sum_{\SSS_n(\g)}\diam\phi(s)\ge \mbox{dist}(\phi(E),\phi(F))$$ et 
$$V_{Q,\rho_n}(X)\lesssim \left(\sup_{\SSS(\EEE_n)}\diam s\right)^{Q-2}\cdot \sum_{s\in\SSS(\EEE_n)} (\diam s)^2
\lesssim\left(\sup_{\EEE_n}\diam s\right)^{Q-2}\mbox{Aire}(\cbar) $$
qui tend vers $0$ d\`es que $Q>2$.\endp

\begin{rema} Ce dernier argument montre qu'en fait, si $X$ est $Q$-Ahlfors r\'egulier pour un $Q>2$, et 
si la suite de quasi-empilements est
conforme, alors le $Q$-module de tout condensateur est nul, et les $Q$-modules combinatoires associ\'es tendent vers $0$ 
par la Proposition \ref{disoka}. Une condition n\'ecessaire pour avoir la conformit\'e d'une suite est donc
que le $Q$-module combinatoire de tout condensateur tend vers $0$ pour tout $Q>2$.\end{rema}

Modulo la construction des quasi-empilements par tuiles de valence born\'ee dont le nerf contient une triangulation
\'equivalente ({\it cf.} \S4 et \S6 dans \cite{bk1} et la Remarque \ref{rmq:tuilevsboules}), on r\'ecup\`ere ainsi des r\'esultats de M.\,Bonk et B.\,Kleiner
 (\Th 1.1 et 1.2 de \cite{bk1})\,:
 
\begin{theo}[M.\,Bonk et B.\,Kleiner]\label{cbk} Si $X$ est hom\'eomorphe \`a
$\SS^2$, lin\'eairement localement connexe et $Q$-r\'egulier, et si on suppose que, ou bien $Q=2$ ou 
bien $X$ est $Q$-loewnesque, alors $X$ est quasisym\'etrique
\`a $\cbar$ (et $Q=2$).\end{theo}

\subsection{Remarque sur l'utilisation des toitures}

Tous nos r\'esultats sont \'etablis pour des toitures. En fait, cette condition n'est essentiellement
utilis\'ee que pour comparer les modules $\hmod_{*}$ et $\mod_*$ 
({\it cf.} (\ref{eqn:trans})) et pour obtenir que les
suites sont conformes. Si on s'affranchit des toitures,
tous les r\'esultats restent vrais en rempla\c{c}ant $\mod_*$ par $\hmod_*$ dans les hypoth\`eses
et les conclusions. Seule la d\'emonstration du \Th \ref{bk10.1} requiert quelques changements car $\phi_n^{-1}$
ne pr\'eserve plus les continua {\it a priori}\,; 
une m\'ethode est de se ramener au Corollaire \ref{th11.1}\,: 
on consid\`ere  $x\in X$ et $r< R\le\diam X$ et on veut montrer 
qu'il existe une constante $C$ qui ne d\'epend que de  $R/r$
telle que $\hmod_2(B(x,r),X\setminus B(x,R),\SSS_n)\le C$.  
On note $E=B(x,r)$ et $F=X\setminus B(x,R)$\,; si $R/r$ est assez
grand, alors la connexit\'e locale lin\'eaire de $X$ loewnesque permet de supposer $E$ et $F$ connexes. 
En partant de (\ref{eqn:1}), 
on suppose que $\De(\phi_n(E),\phi_n(F))\le \de$,
et l'on veut montrer que $\de$ ne peut \^etre trop petit en construisant de nouveaux condensateurs qui feront intervenir la condition de Loewner
pour minorer $\de$. 
On obtiendra ainsi un bon contr\^ole du $2$-module
en utilisant que $\cbar$ est $2$-r\'egulier.

\subsection{Sur les d\'emonstrations des diff\'erentes versions du th\'eor\`eme de Riemann}

Dans tous ces r\'esultats, outre l'hypoth\`ese topologique d'\^etre une sph\`ere, 
on part d'hypoth\`eses qui imposent la pr\'esence de nombreuses courbes sur la surface.
Ces courbes sont quantifi\'ees gr\^ace \`a la notion de modules, et ceux-ci sont utilis\'es pour en d\'eduire des propri\'et\'es
g\'eom\'etriques de la surface.
On suppose donc toujours une
condition g\'eom\'etrico-combinatoire extr\^eme\,: dans les hypoth\`eses de s\'eparation ou 
de conformit\'e, 
on se donne la combinatoire de la sph\`ere de Riemann, dans l'hypoth\`ese $2$-Ahlfors r\'eguli\`ere, 
on se donne
une borne sup\'erieure sur la g\'eom\'etrie des modules, et sous la condition  loewnesque, 
on obtient une borne inf\'erieure. 
De plus, dans les deux approches,
les auteurs imposent des estimations uniformes qui impliquent que les structures complexes
qu'ils d\'efinissent sont \`a distance \qc born\'ee les unes des autres.

\medskip

L'approche de J.W.\,Cannon {\it et al.} est combinatoire et locale,
celle de M.\,Bonk et B.\,Kleiner analytique et globale.
Les d\'emonstrations que nous avons pr\'esent\'ees s'inspirent des deux \'ecoles
et permettent de distinguer deux probl\`emes diff\'erents\,: le premier est de munir la
surface consid\'er\'ee d'une structure complexe compatible avec les modules combinatoires, 
le second consiste \`a montrer que la structure m\'etrique ainsi obtenue est
quasisym\'etrique \`a l'originale 
(lorsque celle-ci est fix\'ee). 
Le second probl\`eme s'interp\`ete alors 
comme une quantification/habillage d'une approche plus combinatoire.

Nous sommes donc partis de l'existence d'une suite de recouvrements dont les
nerfs sont essentiellement des triangulations et dont l'existence est assur\'ee
par M.\,Bonk et B.\,Kleiner sous de bonnes conditions (bien que nous n'ayons pas
besoin de toutes leurs propri\'et\'es, il n'est pas clair comment simplifier
leur construction). 

L'objet du \Th \ref{comp} dont les hypoth\`eses de s\'eparation sont inspir\'ees
des travaux de J.W.\,Cannon {\it et al.} \cite{cfp2}
est de r\'epondre au premier
probl\`eme par un \th de Riemann combinatoire \`a la Cannon. 
Pour cela, on s'appuie sur la th\'eorie des empilements de cercles,
introduite dans ce contexte par  M.\,Bonk et B.\,Kleiner, qui nous ram\`ene 
rapidement sur la sph\`ere de Riemann. Pour ces empilements de cercles, 
la consid\'eration des $\mod_{\inf}$ et $\mod_{\sup}$ de J.W.\,Cannon
et d'arguments assez classiques d'estimations de modules de courbes permet de 
 montrer simplement que les modules combinatoires et analytiques sont 
\'equivalents (Corollaire \ref{inegtr}).  L'int\'er\^et de ces modules est de confronter
deux familles de courbes transverses dont les modules co\"{\i}ncident 
dans le cas analytique (voir aussi le Lemme \ref{infsup}).

On peut alors conclure par le \th d'Ascoli pour obtenir 
une uniformisation {\it globale} par un \homeo
limite qui \'etablit aussi la conformit\'e de sous-suites. 

Pour voir que l'\homeo limite est bien quasisym\'etrique, et r\'epondre au second
probl\`eme, on joue sur les propri\'et\'es
des recouvrements qui d\'ecoulent de leur forme et de la g\'eom\'etrie de l'espace \'etudi\'e\,:
en particulier, on profite que la sph\`ere de Riemann est \`a la fois $2$-r\'eguli\`ere
et $2$-loewnesque pour avoir les in\'egalit\'es compl\'ementaires \`a celles admises
sur $X$ ({\it cf.} Prop.\,\ref{regqs}).

\medskip

Il semblerait que les deux approches originales aient plus
de points communs qu'il n'y para\^{\i}t. 
Par exemple, J.W.\,Cannon consid\`ere les 
modules dans la surface proprement dite et, uniquement
en dimension deux, alors
qu'ils sont d\'efinis sur le nerf par M.\,Bonk et B.\,Kleiner, et en dimension $Q$ variable\,; 
en pratique, ils  se trouvent \^etre \'equivalents dans les m\^emes dimensions 
(Remarque \ref{rmq:tuilevsboules} et Proposition \ref{tvb}). 

De plus, J.W.\,Cannon {\it et al.} \'etablissent dans \cite{ca,cfp2} des
 in\'egalit\'es de type  ``2-r\'eguli\`ere''
 qui tiennent un r\^ole pr\'epond\'erant dans leur argumentation, et 
se ram\`enent ainsi au Corollaire \ref{bk10.4}. 
Se donner des estimations \`a la fois sur les $\mod_{\sup}$ et les $\mod_{\inf}$ 
revient \`a se donner \`a la fois des majorations et des minorations, 
donc des propri\'et\'es de type loewnesque
et de $2$-r\'egularit\'e. 
D'un autre c\^ot\'e, lorsque M.\,Bonk et B.\,Kleiner supposent que $X$ est
$2$-r\'egulier ou se donnent les hypoth\`eses du Corollaire \ref{th11.1}, 
ils supposent en fait que les approximations
sont s\'eparantes, et se rapprochent ainsi du \th de Riemann combinatoire. 

\medskip

Terminons par une question dont la r\'eponse ne semble pas d\'ecouler de la pr\'esente analyse\,:

\begin{ques} Un espace m\'etrique $Q$-r\'egulier $X$  hom\'eomorphe \`a $\SS^2$, $Q>2$, et muni d'une suite
de recouvrements par pi\`eces uniform\'ement rondes, de valence born\'ee, dont le nerf contient une triangulation \'equivalente,
 et conforme, est-il quasisym\'etrique \`a $\cbar$ ?\end{ques}

\section{Applications \`a la conjecture de Cannon}

Ces probl\`emes de discr\'etisation que nous avons trait\'es jusqu'ici sont motiv\'es par la th\'eorie g\'eom\'etrique
des groupes et la g\'eom\'etrisation des vari\'et\'es de dimension $3$ \cite{cap}.
En g\'en\'eral, la classe de quasi-isom\'etrie d'un groupe hyperbolique est d\'etermin\'ee par la topologie
et la structure quasiconforme de son bord. En petite dimension, on tend \`a penser que la topologie
seule suffit \`a d\'eterminer le groupe \`a quasi-isom\'etrie pr\`es.
D'autre part, la conjecture de g\'eom\'etrisation W.P.\,Thurston implique que le groupe fondamental
d'une vari\'et\'e compacte sans bord est hyperbolique si et seulement si la vari\'et\'e peut \^etre munie d'une structure
hyperbolique compl\`ete \cite{th}.

 A l'intersection de ces deux probl\'ematiques, J.W.\,Cannon propose la conjecture suivante dans \cite{cap}.

\begin{conj}[de Cannon] Un groupe hyperbolique de bord hom\'eomorphe \`a $\SS^2$ op\`ere g\'eom\'etriquement sur l'espace
hyperbolique $\HH^3$.\end{conj}

Un \th de  D.\,Sullivan permet de reformuler cette conjecture sous une forme plus analytique.

\begin{conj} Si le bord d'un groupe hyperbolique est hom\'eomorphe \`a $\SS^2$, alors
il est quasisym\'etrique \`a $\cbar$.\end{conj}

\subsection{Approche de Cannon {\it et al.}}

Soit $G$ un groupe hyperbolique au sens de Gromov dont le bord est hom\'eomorphe \`a la sph\`ere $\SS^2$.
On se fixe un syst\`eme fini de g\'en\'erateurs et on lui associe son graphe de Cayley $(X,o)$
vu comme espace m\'etrique. 
On munit son bord $\partial X$ d'une m\'etrique visuelle de param\`etre $\ep$.
Rappelons que pour un rayon  $R>0$ fix\'e assez grand et $x\in X$,
l'ombre $\mho_o(x,R)\subset\partial X$ de la boule $B(x,R)$
port\'ee par une source de lumi\`ere en $o$ ressemble \`a une boule de rayon $e^{-\ep |x|}$.
Pour chaque $n\in\N\setminus\{0\}$, l'ensemble $\mho_n=\{\mho_o(x,R), |x|=n\}$ d\'etermine un 
quasi-empilement de $\partial X$ de valence born\'ee.

J.W.\,Cannon et E.L.\,Swenson ont montr\'e qu'un groupe v\'erifiait la conjecture de Cannon si et seulement
si la famille des recouvrements $\{\mho_n\}_n$ est conforme \cite{cs}. 
Ult\'erieurement, J.W.\,Cannon, W.\,Floyd et W.\,Parry ont affaibli la condition sur $\{\mho_n\}$ 
d'\^etre conforme pour ne supposer que
la variante suivante de la s\'eparation\,: 
pour chaque point $x$ du bord et chaque voisinage $V$ de $x$, il existe un anneau $A$ qui s\'epare $x$
du compl\'ementaire de $V$ et une constante $m=m(A)>0$ tels que $\mod_{\inf}(A,\mho_n)\ge m$ (\Th 8.2 de \cite{cfp2}).  

Dans leur approche, ils utilisent qu'on a {\it a priori} une structure combinatoire \`a l'infini qui est 
``la bonne'' si jamais le groupe op\`ere g\'eom\'etriquement sur $\HH^3$. Du coup, nul besoin d'\'etudier les
propri\'et\'es de la jauge conforme \`a l'infini.

\subsection{Approche de Bonk et Kleiner}

Dans \cite{bk0,bk2},  M.\,Bonk et B.\,Kleiner travaillent avec des hypoth\`eses de nature analytique.
Tout d'abord, il faut savoir que le bord d'un groupe hyperbolique $G$ est toujours Ahlfors-r\'egulier et  lin\'eairement
localement connexe lorsqu'il est connexe. De plus, l'action des isom\'etries 
est uniform\'ement quasim\"obius sur le bord \`a l'infini (les birapports de quatre points sont quasipr\'eserv\'es). 

M.\,Bonk et B.\,Kleiner montrent que s'il existe une m\'etrique
Ahlfors-r\'eguli\`ere de dimension minimale $Q\ge 2$ (parmi les m\'etriques Ahlfors-r\'eguli\`eres de la jauge conforme) 
sur $\partial G$, alors la conjecture est satisfaite.
Cette condition permet de montrer qu'un espace tangent (faible) $T$ de $\partial G$
admet une famille de courbes de $Q$-module non nul \cite{kl}. Or $T$ est aussi $Q$-Ahlfors-r\'egulier et quasim\"obius \'equivalent
\`a un \'epointement de $\partial G$. Ceci implique que $\partial G$ admet aussi une famille de courbes de $Q$-module non nul
(quasi-invariance du module par transformation quasim\"obius \cite{ty}).

La dynamique par transformations quasim\"obius permet d'exploiter cette famille de courbes pour montrer 
que $\partial G$ v\'erifie la condition de Loewner restreinte aux paires de boules
de m\^eme rayon.  M.\,Bonk et B.\,Kleiner montrent alors que cette condition est suffisante pour montrer que $\partial G$
est loewnesque,
et ils en d\'eduisent donc la conjecture en appliquant le Corollaire \ref{cbk}.

\subsection{En guise de conclusion}

Dans les deux articles \cite{cfp2} et \cite{bk2}, les auteurs partent d'une condition sur les modules de courbes\,: 
il existe une/des famille(s) de courbes de module (combinatoire) non nul\,: cette hypoth\`ese -- tr\`es forte --
impose la pr\'esence d'une quantit\'e importante de courbes. Ensuite, la dynamique est utilis\'ee pour propager
cette information \`a tous les endroits et \`a toutes les \'echelles. Cette donn\'ee se traduit par l'existence
d'une structure analytique sur $\partial G$ compatible avec ces modules et l'action de $G$. Enfin, le \th
de Sullivan produit l'action par transformations de M\"obius.

On peut tirer la le\c{c}on suivante\,: la dynamique (cocompacte et \`a distorsion born\'ee)
permet d'obtenir des estimations uniformes \`a toutes les \'echelles et dans toutes les positions \`a partir 
d'estimations plus qualitatives.

\medskip

Il semblerait que l'analogie s'arr\^ete ici. 

Si l'on part du fait qu'il existe une famille de courbes de $Q$-module non nul dans un espace $Q$-r\'egulier ($Q\ge 2$), 
alors l'argument dynamique 
de  M.\,Bonk et B.\,Kleiner \'evoqu\'e plus haut permet de ``balader'' cette famille de courbes un peu partout. De plus,
le lemme du pont et la monotonie des modules combinatoires en fonction de l'exposant $Q$ 
montrent que les $2$-modules de condensateurs sont minor\'es, mais cela n'implique pas {\it a priori} la s\'eparation combinatoire
(il s'agit de l'in\'egalit\'e inverse), ou alors, il faudrait montrer que ces estimations s'appliquent aux familles
de courbes qui s\'eparent les condensateurs ({\it cf.} $\mod_{\inf}$). 

R\'eciproquement, si les hypoth\`eses de s\'eparation combinatoire sont v\'erifi\'ees, on n'a aucune raison d'avoir l'existence
d'une famille de courbes de module non nul\,: il se peut m\^eme qu'il n'y ait aucune courbe rectifiable sur la surface !
N\'eanmoins, si on se donne les hypoth\`eses de s\'eparation combinatoire sur $\partial G$ pour un exposant $Q\ge 2$ et que
l'on suppose aussi que $\partial G$ est $Q$-Ahlfors r\'egulier, alors la Proposition \ref{disoka} montre qu'il existe une
famille de courbes de module positif\,: on est donc ramen\'e \`a \cite{bk2} ! Mais les hypoth\`eses sont tr\`es contraignantes
puisque le $Q$-module combinatoire d'une famille de courbes est toujours plus petit que son $2$-module, et tend
m\^eme vers $0$ dans les cas favorables !

\medskip

En r\'esum\'e, une premi\`ere impression indique que l'approche de M.\,Bonk et B.\,Kleiner est plus
approfondie puisque l'on r\'ecup\`ere des propri\'et\'es combinatoires (construction de triangulations et contr\^ole
des modules combinatoires associ\'es)
\`a partir de donn\'ees analytiques, alors
que l'approche de J.W.\,Cannon {\it et al.} part directement de donn\'ees combinatoires. En revanche, il est
plus raisonnable de consid\'erer des modules combinatoires, puisque l'on ne sait pas en g\'en\'eral 
s'il existe des courbes rectifiables...
De plus, il faut bien
comprendre que la donn\'ee initiale est celle d'un groupe hyperbolique de bord $\SS^2$, donc les donn\'ees
sont plut\^ot de nature combinatoire et topologique, et les propri\'et\'es analytiques ne pourront
provenir  -- comme c'est d\'ej\`a le cas -- que de consid\'erations combinatoires... 

\medskip

Enfin, notons que  dans les deux cas, les conditions suffisantes mises
en avant (existence quantifi\'ee de courbes) sont {\it a priori} et en g\'en\'eral difficiles -- 
pour ne pas dire quasi-impossibles -- \`a v\'erifier...
ce qui laisse le probl\`eme encore ouvert ! 

En effet, on sait qu'il existe des m\'etriques $Q$-r\'eguli\`eres sur $\partial G$, pour des exposants $Q> 2$\,: on
sait donc  comment majorer les $Q$-modules de condensateurs, 
mais comment montrer qu'on a une famille de $Q$-module positif ? ou que sa dimension est minimale ?? De m\^eme,
qu'est-ce qui nous permet de contr\^oler les $2$-modules combinatoires ? Comment montrer que l'on peut
choisir $Q=2$ sans r\'epondre aux questions pr\'ec\'edentes ?

\appendix
\section{Notions de th\'eorie g\'eom\'etrique des fonctions}\label{tgf}

Des pr\'ecisions sur les notions d\'evelopp\'ees dans  cet appendice peuvent \^etre trouv\'ees par exemple dans les ouvrages
\cite{ah2,ah1,he}.

Nous nous pla\c{c}ons dans un premier temps sur la sph\`ere de Riemann. Le module est un outil
fondamental de l'analyse complexe issu de la m\'ethode {\sl longueur-aire}, et ce, pour au moins
deux raisons. Tout d'abord il se traduit souvent par des contraintes g\'eom\'etriques, et ensuite, il se comporte
bien par rev\^etements (ramifi\'es) holomorphes (ou quasiconformes).

\begin{defi}[module (analytique) de familles de courbes]\label{defi:modana}
Soit $\Gamma$ une famille de courbes de $\cbar$ et $Q > 1$ un
r\'eel. On d\'efinit le $Q$-module de $\Gamma$ par 
$$
\mod_{Q}\Gamma = \inf \int_{\cbar} \rho^{Q} 
$$
o\`u l'infimum est pris sur toutes les fonctions bor\'eliennes (dites admissibles)
$\rho : X \rightarrow [0, \infty]$ telles que, pour toute courbe
rectifiable $\gamma \in \Gamma$, $\int_{\gamma} \rho ds \geq 1$.\end{defi}

Lorsque $X=\cbar$ et $Q=2$, on obtient ainsi un invariant conforme.

Donnons quelques propri\'et\'es \'el\'ementaires du module.
\begin{enumerate}
\item $\mod_{Q} (\emptyset) = 0$;
\item si $\Gamma_{1} \subset \Gamma_{2}$, $\mod_{Q} \Gamma_{1} \leq
\mod_{Q} \Gamma_{2}$;
\item $\mod_{Q} ( \bigcup_{i=1}^{\infty} \Gamma_{i}) \leq
\sum_{i=1}^{\infty} \mod_{Q} \Gamma_{i}$ ;
\item si $\Gamma_{1}$ et $\Gamma_{2}$ sont deux familles de courbes
telles que toute courbe $\gamma_{1}$ dans $\Gamma_{1}$  poss\`ede une
sous-courbe $\gamma_{2} \in \Gamma_{2}$, alors $\mod_{Q}
\Gamma_{1} \leq \mod_{Q} \Gamma_{2}$.
\end{enumerate}

On a le crit\`ere suivant.

\begin{prop}[crit\`ere de Beurling] Soient  $\G$ une famille de courbes.
Une m\'etrique $\rho$ est extr\'emale s'il existe une sous-famille
$\G_0\subset\G$ telle que \ben
\item pour tout $\g\in\G_0$, $\ell_\rho(\g) = L_\rho(\G)$\,;
\item si $h:\cbar\to\R$ v\'erifie $\int_{\g} h\ge 0$ pour toute courbe $\g\in\G_0$,
alors on a $\int h \rho ^{Q-1}\ge 0$.\een
De plus, cette m\'etrique est unique \`a normalisation pr\`es.\end{prop}
Quitte \`a jeter une sous-famille de courbes de module z\'ero, il existe 
une m\'etrique extr\'emale, mais il n'est pas clair que la r\'eciproque soit vraie en g\'en\'eral.

\Preuve Soit $\si\in\MMM_Q$ que l'on suppose normalis\'ee pour que $L_{\si}(\G)=L_\rho(\G)$. Du coup, 
$\ell_\si(\g)\ge L_\rho(\G)=\ell_\rho(\g)$ pour toute $\g\in\G_0$ d'apr\`es (1). En posant $h=\si-\rho$, on obtient par (2)
$$\int (\si-\rho) \rho^{Q-1}\ge 0$$
soit $$\int \rho^Q\le \int\si\rho^{Q-1} \le \left(\int \si^Q\right)^{1/Q}\cdot\left(\int \rho^Q\right)^{1-1/Q}$$
par l'in\'egalit\'e de H\"older. Du coup, on a $\mod_Q(\G,\rho)\le\mod_Q(\G,\si)$.
Le cas d'\'egalit\'e se produit lorsque l'in\'egalit\'e de H\"older est une \'egalit\'e, donc si $\si^Q$ et 
$(\rho^{Q-1})^{Q/(Q-1)}$
sont proportionnelles, soit si $\si=\rho$ presque partout d'apr\`es leurs normalisations.\endp

\begin{defi}[module d'anneaux, de quadrilat\`eres et de condensateurs]\label{defi:modanquadcon}
Si $A\subset\cbar$ est un anneau, on lui associe la famille $\G_t$ des courbes qui joignent les composantes de bord, et 
$\G_s$  la famille des courbes qui s\'eparent ces composantes. On a (lemme et d\'efinition)
$$\mod\, A =\mod_2\G_s = \frac{1}{\mod_2 \G_t}\,.$$
Un quadrilat\`ere $Q$ est un domaine de Jordan avec quatre points marqu\'es sur le bord. Ces points d\'ecoupent
le bord en quatre ``c\^ot\'es'' $c_1,c_2,c_3,c_4$ cons\'ecutifs. On peut d\'efinir deux modules, selon que
l'on consid\`ere la famille $\G_{\rm impair}$ des courbes qui relient $c_1$ et $c_3$ ou la famille $\G_{\rm pair}$
 des courbes qui relient $c_2$ et $c_4$. Similairement au cas de l'anneau, on a toujours 
 $$\mod_2\G_{\rm impair} = \frac{1}{\mod_2 \G_{\rm pair}}\,.$$
Un condensateur est form\'e de deux continua $E$ et $F$ disjoints. On note $\G=\G(E,F)$ la famille des courbes
qui relient $E$ et $F$, et on d\'efinit $$\mod_Q(E,F)=\mod_Q \G.$$\end{defi}

\begin{rema} Il faut faire attention \`a ce que le module d'un anneau est {\it l'inverse} du module du condensateur
form\'e des composantes du compl\'ementaire.\end{rema}

Le crit\`ere de Beurling permet de montrer que dans un rectangle, le module est le rapport des longueurs
des c\^ot\'es, et dans le cas d'une couronne $A=\{r < |z|<R\}$, $\mod\, A= (1/2\pi)\log R/r$.
En particulier, si un anneau a un tr\`es grand module, alors au moins l'une des composantes connexes de son
compl\'ementaire a un petit diam\`etre (cela d\'ecoule p.ex. du \th de la couronne de Teichm\"uller
qui affirme qu'un anneau de grand module contient toujours une couronne de module comparable,
voir \cite[Thm 4.8]{ah1}).

Les applications qui sont naturelles dans ce contexte sont les \homeos quasiconformes. Ils ont la particularit\'e
de pouvoir \^etre d\'efinis de bien des mani\`eres dans $\C$. Nous en citons quelques unes.

Soit $f:\C\to\C$ une application et 
soient $r>0$ et $z\in\C$\,; on d\'efinit $$\begin{array}{ll} L_f(z,r)= & \sup\{ |f(z)-f(w)|,\,|z-w|\le r\}\,,\\
\ell_f(z,r)= & \inf\{ |f(z)-f(w)|,\,|z-w|\ge r\}\,.\end{array}$$
Si $\G$ est une famille de courbes, alors $f(\G)$ d\'esigne $\{f(\g),\g\in\G\}$.

\begin{defi}[\homeo \qc du plan] Un \homeo $f:\C\to\C$ est \qc si l'une des conditions \'equivalentes suivantes
est v\'erifi\'ee.
\begin{itemize}
\item[{\bf [QC]}]  Il existe $H<\infty$ telle que pour tout $z\in\C$,
$$\limsup_{r\to 0} \frac{L_f(z,r)}{\ell_f(z,r)}\le H\,.$$
\item [{\bf [G1]}] Il existe $K<\infty$ telle que pour toute famille $\G$
de courbes rectifiables, $$(1/K)\cdot \mod_2\G\le \mod_2 f(\G)\le
K\cdot \mod_2\G\,.$$

\item[{\bf [G2]}] Il existe $K<\infty$ telle que pour tout anneau $A\subset\C$,
$$(1/K)\cdot \mod\,A\le \mod\,f(A)\le K\cdot \mod\,A\,.$$

\item[{\bf [QS]}] Il existe une fonction continue et croissante
$\eta:\R_+\to\R_+$ telle que $\eta(0)=0$, et si $z,w_1$ et $w_2$
sont trois points du plan tels que $|z-w_1|\le t\cdot|z-w_2|$ alors $$|f(z)-f(w_1)|\le \eta(t)\cdot|f(z)-f(w_2)|\,.$$

\item[{\bf [FQS]}] Il existe $H<\infty$ telle que pour tout $z\in\C$ et tout $r>0$, $$L_f(z,r)\le H\cdot \ell_f(z,r)\,.$$
\end{itemize}\end{defi}

Pour utiliser des techniques similaires dans des espaces m\'etriques, il convient de rajouter quelques restrictions
de nature g\'eom\'etrique. Une condition naturelle est celle d'\^etre Ahlfors r\'egulier. Elle est implicite
dans les estimations faites sur $\cbar$.

\begin{defi}[Ahlfors r\'egularit\'e] Si $Q>0$, un espace m\'etrique $X$ est $Q$-Ahlfors r\'egulier s'il existe une mesure 
bor\'elienne $\mu$ telle que, pour tout $R \le\diam X$, et tout $x\in X$ et toute boule ferm\'ee
$B(x,R)$, $\mu(B(x,R))\asymp R^Q$.\end{defi}
Dans ce cas, on peut choisir $\mu=\HHH_Q$ la mesure de Hausdorff $Q$-dimensionnelle.

On remarque qu'un tel espace est toujours {\sl doublant} {\it i.e.}, il existe un naturel $N$
tel que tout ensemble born\'e peut \^etre recouvert par au plus $N$ ensembles de diam\`etre moiti\'e.

\medskip 

Une courbe $\gamma$ dans $X$ est une application continue d'un
intervalle compact $I$ de $\R$ dans $X$. On peut, comme dans les espaces
euclidiens, d\'efinir la longueur de $\gamma$. Si cette longueur
$\ell(\gamma)$ est
finie, on dira que la courbe est rectifiable. Dans ce cas, on peut 
param\'etrer $\gamma$ par la longueur d'arc $\gamma_{s}$ et pour toute
fonction bor\'elienne $\rho : X \rightarrow \R_{+}$, on d\'efinit 
$$
\int_{\gamma} \rho ds = \int_{0}^{\ell(\gamma)} \rho \circ \gamma_{s} (t)
dt.
$$
Cela suffit pour \'etendre la notion de module (analytique) de courbes
aux espaces mesur\'es.

On conclut ces propri\'et\'es par un r\'esultat de S.\,Keith (Proposition 4.2.1 \cite{ke}) qui nous sera utile dans l'appendice suivant.

\begin{prop}[S.\,Keith]\label{modcont} 
Soient $(X,\mu)$ un espace m\'etrique compact de masse finie et  $Q>1$. Pour toute famille de courbes $\G$ et tout
$\ep>0$, il existe une fonction continue $g:X\to\R_+$ telle que $L_g(\G)=1$ et $$\int_X g^Q\le \mod_Q\G +\ep\,.$$
\end{prop}

Le prochain lemme \'etablit une borne sup\'erieure sur le $Q$-module dans un espace $Q$-r\'egulier.

\medskip

\begin{lemm}\label{inegAR} Si $(X,\mu)$ est $Q$-Ahlfors r\'egulier, $x\in X$ et  $\ell <L< \diam X$ avec $L\ge 2\ell$, alors 
$$\mod_Q (B(x,\ell),X\setminus B(x,L))\lesssim (1/\log L/\ell )^{Q-1}$$\end{lemm}

La d\'emonstration consiste \`a estimer le module en utilisant la m\'etrique
$$\rho(y)= (1/\mbox{dist\,} (x,y) )\chi_{B(x,L)\setminus B(x,\ell)}(y)\,.$$ 

L'estimation oppos\'ee du module n'est {\it a priori} pas automatique sous la seule condition de r\'egularit\'e.
Son existence provient de ce que l'on appelle la condition de Loewner, qui a \'et\'e introduite par
J.\,Heinonen et P.\,Koskela \cite{hk1} dans un contexte purement m\'etrique.
Dans ces espaces, les
m\'ethodes de discr\'etisation sont particuli\`erement pertinentes.

\begin{defi}[distance relative]\label{dr} Si $(E,F)$ est un condensateur dans un espace m\'etrique, on d\'efinit
leur distance relative par
$$\Delta (E,F) := \frac{\hbox{\rm dist} (E,F) }{\min \{\diam E, \diam F\}}.$$\end{defi}

\begin{defi}[espace loewnesque] Soit $(X,\mu)$ un espace m\'etrique, mesur\'e de dimension de
Hausdorff $Q$. 
On dit que $X$ est un espace loewnesque s'il v\'erifie la condition de Loewner suivante\,:
il existe
une fonction $\psi : ]0,\infty [ \rightarrow ]0,\infty [$ telle que,
pour tout condensateur $(E,F)$ et
 pour tout $t\geq \Delta (E,F)$, on ait
$$
\mod_{Q} (E,F) \geq \psi (t).
$$\end{defi}

Un espace loewnesque Ahlfors r\'egulier v\'erifie de nombreuses propri\'et\'es similaires aux espaces euclidiens.
Par exemple, 
un tel espace est toujours lin\'eairement localement connexe et quasiconvexe {\it i.e.},

\begin{defi}[connexit\'e locale lin\'eaire] Un espace m\'etrique 
$X$ est lin\'eairement localement connexe s'il existe $C > 0$ telle que, pour tout $x \in X$, tout $r > 0$, on ait\,:
\ben \item tout couple de points dans $B(x,r)$ appartient \`a un continuum contenu dans $B(x,Cr)$;
\item  tout couple de points dans $X \setminus \bar{B} (x,r)$ appartient \`a un continuum contenu dans $X \setminus \bar{B} (x, (1/C) r )$.\een\end{defi}

\begin{defi}[quasiconvexit\'e] Un espace m\'etrique 
$X$ est quasiconvexe s'il existe $C > 0$ telle que, pour tout $x,y \in X$, il existe
une courbe $\g$ qui relie $x$ \`a $y$ telle que $\ell(\g)\le C|x-y|$.\end{defi}

Dans ce contexte, on peut aussi d\'efinir la classes des \homeos \qcs en consid\'erant les applications qui v\'erifient
la condition [QC], la classe des \homeos g\'eom\'etriquement \qcs qui v\'erifient [G1] ou [G2],
la classe quasisym\'etrique, qui v\'erifient
la condition [QS], ou faiblement quasisym\'etriques [FQS]...

D\'efinis sur des espaces doublants et connexes, les \homeos faiblement quasisym\'etriques sont quasisym\'etriques.
Dans un espace de Loewner Ahlfors r\'egulier, toutes ces d\'efinitions sont encore \'equivalentes (ce n'est pas
le cas en g\'en\'eral). La condition loewnesque permet notamment de passer du local au global \cite{hk1}.

\section{Discr\'etisation dans les espaces m\'etriques}\label{ap:dis}

On \'etend la notion de modules combinatoires \`a des espaces m\'etriques plus g\'en\'eraux. 
Dans cet appendice,
on supposera que $(X,\mu)$ est un espace m\'etrique connexe compact localement connexe $Q$-Ahlfors 
r\'egulier pour un
$Q>1$ de diam\`etre $1$. L'objectif de ce paragraphe est de montrer que l'on peut discr\'etiser 
les modules de mani\`ere
raisonnable dans ce contexte (voir Proposition \ref{disoka}).

On se place dans la situation suivante.

\begin{defi}[Quasi-empilement]\label{def:qemp} Si $X$ est un espace m\'etrique, on dira qu'un recouvrement $\SSS$ est un
{\sl quasi-empilement} s'il existe une constante $K\ge 1$ et, pour chaque $s\in \SSS$, un point $x_s$
de $s$ et une taille $r_s>0$  tels que \ben
\item $B(x_s,r_s)\subset s\subset B(x_s,K\cdot r_s)\,.$
\item Chaque boule interne $B(x_s,r_s)$ intersecte au plus $K$ autres boules internes 
$B(x_{s'},r_{s'})$.\een\end{defi}
Un empilement de cercles n'est pas un quasi-empilement, puisqu'il ne recouvre pas l'espace.
On justifie la terminologie par le fait
qu'un quasi-empilement se souvient de propri\'et\'es d'un empilement\,: les
 pi\`eces sont \`a peu pr\`es rondes, et 
(2) est une extension naturelle dans notre contexte que les boules internes soient deux \`a deux disjointes.

\medskip

On se donne une suite de quasi-empilements $(\SSS_n)$ dont la maille tend vers $0$.
On peut proc\'eder ainsi\,:
pour tout $r\in ]0,1]$ on se fixe  une famille de boules ferm\'ees disjointes de 
$\BBB(r)=\{B_i=B(x_i,r/5)\}_i$ de rayon $r/5$ telle que $\hat{\BBB}(r)=\{\hat{B}_i=B(x_i,r)\}$ recouvre $X$. Puisque 
$X$ est Ahlfors r\'egulier, il existe une constante $K\in\N$, ind\'ependante de $r$, telle que chaque boule
de $\hat{\BBB}(r)$ intersecte au plus $K$ autres boules. En effet, si $N$ est le nombre de boules $\hat{B}_j$ qui intersectent
une boule $\hat{B}$ donn\'ee, alors on a $B_j\subset 3\hat{B}$ donc  
$$r^Q\gtrsim \mu(3\hat{B})\ge \sum_{1\le j\le N} \mu(B_j)\gtrsim N r^Q\,.$$

Notons enfin que chaque $\hat{\BBB}(r)$ est un quasi-empilement et on peut prendre 
$\SSS_n=\hat{\BBB}(1/2^n)$ et $\BBB_n= \BBB(1/2^n)$ avec $K=5$.

\begin{prop}\label{disoka} Soit $(X,\mu)$ un espace $Q$-r\'egulier, $Q >1$, muni d'une suite
de quasi-empilements dont la maille tend vers z\'ero.
Pour $L>0$, on note $\G_L$ les courbes de $X$ de diam\`etre au moins $L$.
Pour  $n$ assez grand, on a
$$\mod_Q(\G_L,\SSS_n)\asymp \mod_Q\G_L$$
si $\mod_Q \G_L>0$ et sinon, $\lim \mod_Q(\G_L,\SSS_n)=0$.
De m\^eme, pour tout condensateur $(E,F)$ et pour $n$ assez grand, on a
$$\mod_Q(E,F,\SSS_n)\asymp \mod_Q(E,F)$$
si $\mod_Q(E,F)>0$ et sinon, $\lim \mod_Q(E,F,\SSS_n)=0$.
\end{prop}

Une des in\'egalit\'es sera un corollaire quasi-imm\'ediat du lemme suivant, que l'on
extrait de la proposition car il est plus g\'en\'eral\,: 

\begin{lemm}\label{pongen} Soit $(X,\mu)$ un espace $Q$-r\'egulier muni d'un $K$-quasi-empilement
$\SSS$, et soit $\G$ une famille de courbes. On suppose qu'il existe une constante $\kappa > 0$
telle que, pour tout $s\in\SSS$ et toute $\g\in\G$, si $\g \cap s\ne \emty$ alors
$\diam ( \g\cap (2K)\cdot B(s))\ge \kappa\cdot \diam B(s)$,
o\`u $B(s)$ est la boule interne de $s$.

Alors, on a $$\mod_Q\G\lesssim \mod_Q(\G,\SSS)\,.$$
\end{lemm}

Avant de d\'emontrer ces \'enonc\'es, on rappelle un lemme qui nous sera utile.

\REFLEM{lmclan} Soit $(X,\mu)$ un espace m\'etrique mesur\'e qui v\'erifie la condition de doublement
de volume. On consid\`ere une famille de boules ${B\in\BBB}$ et on associe \`a chacune un poids $a_B>0$.
Pour tout $p>1$, pour tout $\la\in ]0,1[$, il existe une constante $C=C(p,\la)>0$ ind\'ependante
de $\BBB$ et des poids, telle que 
$$\int \left(\sum a_B\chi_B\right)^p d\mu\le C \int \left(\sum a_B\chi_{\la B}\right)^p d\mu\,.$$
\ENDLEM

La d\'emonstration de ce lemme utilise les fonctions maximales de Hardy-Littlewood et consiste \`a 
estimer la norme de cette fonction en tant qu'\'el\'ement du dual de $L^q$, o\`u $(1/p)+(1/q)=1$ (voir
\cite{bo} dans le cadre euclidien).

\medskip

{\noindent\sc D\'emonstration du Lemme \ref{pongen}.} 
Si $\rho\in\MMM_Q(X,\SSS)$, on d\'efinit la soupe $\hat{\rho}$ du poids $\rho$ par
$$\hat{\rho}=\sum_{s\in\SSS}\frac{\rho(s)}{\diam B(s)} \chi_{2KB(s)}\,,$$
o\`u $\chi_{2KB(s)}$ d\'esigne la fonction caract\'eristique de $2KB(s)$.
Du coup, si $\g\in\G$, alors
$$\begin{array}{ll}\ell_{\hat{\rho}}(\g) & \ge \dis\sum_{s\in\SSS(\g)} \dis\int_{\g\cap 2KB(s)} \frac{\rho(s)}{\diam s}\\&\\
& \ge \kappa\dis\sum_{s\in\SSS(\g)} \rho(s)\\ &\\
& \ge \kappa L_\rho(\G,\SSS)\,.\end{array}$$

\medskip

D'autre part,
$$\begin{array}{ll} \mbox{V}_Q(\hat{\rho}) & = \dis\int_{X} \left(\dis\sum_{s\in\SSS} \frac{\rho(s)}{\diam s}\chi_{2KB(s)}\right)^Q\\&\\
& \lesssim\dis\int_{X} \left(\dis\sum_{s\in\SSS} \frac{\rho(s)}{\diam s}\chi_{B(s)}\right)^Q\end{array}$$
par le Lemme \ref{lmclan}. 
Or, pour tout $x\in X$, on a
$$\begin{array}{ll}\left(\dis\sum_{s\in\SSS} \dis\frac{\rho(s)}{\diam s}\chi_{B(s)}(x)\right)^Q  &
\le  \left( (K+1)\cdot \max \left\{\dis\frac{\rho(s)}{\diam s}\chi_{B(s)}(x)\right\}\right)^Q \\&\\
& \le (K+1)^Q\dis\sum_{s\in\SSS}  \left( \dis\frac{\rho(s)}{\diam s}\chi_{B(s)}(x)\right)^Q\,.\end{array}$$
Par cons\'equent, 
on a
$$\begin{array}{ll} \mbox{V}_Q(\hat{\rho})
& \lesssim\dis\sum_{s\in\SSS}\dis\int_{B(s)} \left( \dis\frac{\rho(s)}{\diam s}\right)^Q\\&\\
& \lesssim\dis\sum_{s\in\SSS}\rho(s)^Q\end{array}$$ car $\mu(B(s))\asymp \diam B(s)^Q$. Donc
$$V_Q(\hat{\rho})\lesssim V_Q(\rho,\SSS)$$  et on obtient ainsi $$\mod_Q\G\lesssim\mod_Q(\G,\SSS)\,.$$ 
\endp

Pour l'autre in\'egalit\'e, nous allons  adapter l'argument de la Proposition 3.2.1 de \cite{kl}.
Dans ce but, on d\'efinit {\it une} $\SSS_n$-approximation d'une courbe $\g:[0,1]\to X$
comme une suite finie de pi\`eces $\{s_1,\ldots,s_k\}$ de $\SSS_n$ deux \`a deux distinctes, telle qu'il existe 
$0=t_0<t_1<\ldots <t_j<\ldots <t_k=1$ qui v\'erifient $\{\g(t_j),\g(t_{j+1})\}\subset s_{j+1}$. 
On note $\ell_n(\g)=\sum_{0\le j< k} d(\g(t_j),\g(t_{j+1}))$ et on d\'efinit l'application 
$a_n(\g):[0,\ell_n(\g)]\to X$, constante par morceaux, par
$$a_n(\g)\left(\left[\sum_{0\le j< m} d(\g(t_j),\g(t_{j+1})), \sum_{0\le j\le m}d(\g(t_j),\g(t_{j+1})) \right[\right)= \g(t_{m})$$
pour $0\le m \le k-1$ et $a_n(\g)(\ell_n(\g))=\g(1)$.

Notons qu'en g\'en\'eral il existe plusieurs $\SSS_n$-approximations d'une m\^eme courbe. Cependant, on a toujours
$\ell_n(\g)\le \sum_{\SSS_n(\g)}\diam s$.

\medskip

{\noindent\sc D\'emonstration de la Proposition \ref{disoka}.} 
On se fixe un condensateur $(E,F)$ et un diam\`etre $L>0$, et on note $\G=\G(E,F)$ ou $\G=\G_L$. 

On montre dans un premier temps que $\mod_Q\G\lesssim\mod_Q(\G,\SSS_n)$ 
pour $n$ assez grand, o\`u les constantes implicites ne d\'ependent ni de $n$, ni du condensateur.
Par hypoth\`eses, il existe une constante $K\ge 1$
et  pour chaque $s\in\cup\SSS_n$
une boule $B(s)$ telles que $B(s)\subset s\subset K\cdot B(s)$. Comme la maille de $\SSS_n$ tend vers $0$,
il existe $\kappa>0$ tel que, pour tout $n$ assez grand, pour toute $\g\in\G$ et toute  $s\in\SSS(\g)$, 
$$\diam ( \g\cap (2K)\cdot B(s))\ge \kappa\cdot \diam B(s).$$
Pour ces $n$-ci, on obtient  $\mod_Q\G\lesssim\mod_Q(\G,\SSS_n)$ par le Lemme \ref{pongen}.

\medskip

Quant \`a la r\'eciproque, on se fixe $\ep >0$. 
D'apr\`es la Proposition \ref{modcont},  il existe une fonction continue $g:X\to\R_+$ 
telle que $L_g(\G)=1$ et $$\int_X g^Qd\mu\le\mod_Q\G +\ep\,.$$
On peut supposer $g(x)\ge \de>0$ puisque $X$ est compact et $\mu(X)<\infty$. On d\'efinit, pour $x\in X$, 
$$g_n(x)=\inf\{ g(y),\ y\in s,\ x\in s, \ s\in\SSS_n\}\,.$$
Puisque $g$ est continue et la maille de $\SSS_n$ tend vers $0$, la suite $(g_n)$ tend uniform\'ement vers $g$.

Notons $\rho_n(s):= (3/2)\inf_s g\cdot\diam s$, et 
montrons
que $L_{\rho_n}(\G,\SSS_n)\ge 1$ pour $n$ assez grand.
Si ce n'est pas le cas,
alors il existe une sous-suite $(n_p)$ ainsi que des courbes 
$(\g_p)$ de $\G$ telles que $\ell_{\rho_{n_p}}(\g_p,\SSS_{n_p})<1$. On se fixe des approximations $a_{n_p}(\g_p)$. 
Ceci implique  que 
d'une part $$\sum_{\SSS_{n_p}(\g_p)} \diam s\cdot\inf_{s}g < 1 -(\de/2) \sum_{\SSS_{n_p}(\g_p)}\diam s\,,$$
et d'autre part $\ell_{n_p}(\g_p) <1/\de$.

Notons $\ell=\sup_p \ell_{n_p}(\g_p)$ et \'etendons $a_{n_p}(\g_p)$ \`a $[0,\ell]$ en posant 
$a_{n_p}(\g_p)([\ell_{n_p}(\g_p),\ell])=a_{n_p}(\g_p)(\ell_{n_p}(\g_p))$.

On obtient ainsi une suite relativement compacte pour la convergence uniforme 
dont toute limite est une application $1$-lipschitzienne (argument similaire au \th d'Arz\'ela-Ascoli). Soit
$\g$ l'une d'elles\,: il s'agit, par construction, d'une courbe $\g\in\G$\,: dans le premier cas, elle relie
$E$ et $F$; dans le second, son diam\`etre est bien minor\'e par $L$.

On obtient, par le lemme de Fatou,
$$\begin{array}{ll}\dis\int_{\g} gds & \le \liminf \dis\int_0^L g_{n_p}\circ a_{n_p}(\g_p)\\&\\
& \le \liminf \dis\sum_{\SSS_{n_p}(\g_p)} \diam s\inf_{s}g \\&\\
& \le 1 - (\de/2)\limsup \dis\sum_{\SSS_{n_p}(\g_p)}\diam s \\ & \\
& \le 1 -(\de/2)\ell(\g)\,.\end{array}$$ Ceci contredit $L_g(\G)=1$.

Du coup, pour $n$ assez grand, on a 
$$\begin{array}{ll}\mod_Q(\G,\SSS_n) & \lesssim\dis\sum (\inf_s g)^Q\cdot(\diam s)^Q\\ &\\
& \lesssim  \dis\sum (\inf_s g)^Q\cdot\mu(B(s))\\&\\
& \le \dis\sum \dis\int_{B(s)} g^Qd\mu\\&\\
& \lesssim \dis\int_X g^Qd\mu\\&\\
& \le \mod_Q\G+\ep.\end{array}$$
Par cons\'equent, si $\mod_Q\G=0$, on obtient $\mod_Q(\G,\SSS_n)\to 0$ avec $n$, et si
$\mod_Q\G>0$, on obtient $\mod_Q(\G,\SSS_n) \lesssim \mod_Q\G$ en prenant {\it e.g.} $\ep= \mod_Q\G$.\endp

Cette proposition permet de donner une nouvelle d\'emonstration qu'une application quasim\"obius entre
$Q$-espaces r\'eguliers quasipr\'eserve les $Q$-modules
de condensateurs quantitativement \cite{ty}. Rappelons, suivant J.\,V\"ais\"al\"a \cite{vam},
qu'un \homeo $\phi:X\to Y$ entre espaces m\'etriques
est quasim\"obius s'il existe un \homeo $\eta:\R_+\to\R_+$ tel que, quelque soit le quadruplet
de points distincts $x_1,x_2,x_3,x_4\in X$, on a
$$\frac{|\phi(x_1)-\phi(x_2)|}{|\phi(x_1)-\phi(x_3)|}\frac{|\phi(x_3)-\phi(x_4)|}{|\phi(x_3)-\phi(x_4)|}\le
\eta\left(\frac{|x_1-x_2|}{|x_1-x_3|}\frac{|x_3-x_4|}{|x_3-x_4|} \right)\,.$$

\begin{theo}[J.\,Tyson] Soient $X,Y$ deux espaces compacts $Q$-r\'eguliers et $\phi:X\to Y$ un \homeo quasim\"obius. Alors, pour tout condensateur
$(E,F)$, on a $$\mod_Q(E,F)\asymp \mod_Q(\phi E,\phi F)\,,$$
o\`u les constantes implicites ne d\'ependent que des constantes des propri\'et\'es \'evoqu\'ees sur $X$ et $Y$.\end{theo}

\Preuve On se fixe une suite de recouvrements $(\SSS_n)$ comme pour la Proposition \ref{disoka}. 
Une application quasim\"obius est localement quasisym\'etrique (quantitativement)
donc l'image des $\SSS_n$ reste une suite de recouvrements de valence born\'ee par pi\`eces uniform\'ement rondes dont la maille
tend vers $0$, 
donc, en utilisant deux fois la Proposition \ref{disoka} pour $n$ assez grand, on obtient
$$\mod_Q(E,F)\asymp\mod_Q(E,F,\SSS_n)= \mod_Q(\phi E,\phi F,\phi \SSS_n)\asymp \mod_Q(\phi E,\phi F)$$
ce qui \'etablit le th\'eor\`eme. \endp

\section{Convergence et topologie de Hausdorff}\label{apphaus}

Cet appendice pr\'ecise les notions de convergence utilis\'ees notamment dans la d\'emonstration du
\Th \ref{comp}.

Soit $(Z,d)$ un espace m\'etrique. Si $X,Y\subset Z$, on note $$\partial(X,Y)=\sup_{x\in X} d(x,Y)\,.$$

On a $\partial (X,Y)=0$ si et seulement si $\overline{X}\subset \overline{Y}$.

On d\'efinit  par
$$d_H(X,Y)=\max\{\partial (X,Y),\partial(Y,X)\}$$
la {\sl distance de Hausdorff} entre $X$ et $Y$. On dira qu'une suite $(X_n)_n$ de parties de  $Z$ tend vers $X$ 
si $d_H(X_n,X)\to 0$.

En notant $V_\ep(X)=\{x\in Z,\ d(x,X)<\ep\}$, cela signifie que, pour tout $\ep>0$, il existe $n_0$ tel que si $n\ge n_0$ alors
$X_n\subset V_\ep(X)$ et $X\subset V_\ep(X_n)$.

Rappelons que si $X$ un espace m\'etrique propre, alors l'ensemble des compacts non vides de $X$ muni de la distance
de Hausdorff est un espace m\'etrique complet.

Dans la d\'emonstration du \Th \ref{comp}, les ensembles $\PPP_n$ tendent dans la topologie de Hausdorff vers
$X$. 

\medskip

Dans un premier temps, on r\'einterpr\`ete la convergence au sens de Hausdorff \`a l'aide de la convergence
uniforme de fonctions. Ensuite, on \'etablit un \th de convergence de fonctions d\'efinies sur des sous-espaces
convergents. On se place dans un contexte tr\`es proche de ce dont on a besoin. Cette approche est inspir\'ee
par les travaux de M.\,Gromov sur la m\'etrique de \og Hausdorff-Gromov \fg\ dans \cite{gr}, et nos r\'esultats s'\'etendent \`a
ce contexte plus g\'en\'eral sans probl\`eme particulier.

\begin{lemm}  Si $Z$ est un espace m\'etrique compact, alors, pour tout $\ep>0$, il existe un nombre $N(\ep)$ tel que 
tout compact $K$ de $Z$ peut \^etre recouvert par $N(\ep)$ boules de rayon $\ep$ centr\'ees sur $K$. \end{lemm}

\Preuve Si $\ep>0$ est fix\'e, on consid\`ere un recouvrement fini de $Z$ par des boules de rayon $\ep/2$. On note $N(\ep)$
le nombre de ces boules. Pour chacune
de ces boules $B$, pour chaque $x\in B$, on a $B\subset B(x,\ep)$. Donc, si $K$ est compact,
alors il est recouvert par au plus $N(\ep)$ boules de rayon $\ep$.\endp

\bigskip

Soit $Z$ un espace m\'etrique compact. On consid\`ere la suite $(N_k)_{k\ge 1}$ des entiers obtenus par le lemme pr\'ec\'edent en prenant
$\ep=1/2^k$. 

On note $$A_k=\left\{(n_j)_{1\le j\le k}\in\prod_{1\le j\le k} [1,N_j]\right\}\quad\mbox{ et } p_k:A_{k+1}\to A_k$$ l'application 
$(n_j)_{1\le j\le k+1}\mapsto  (n_j)_{1\le j\le k}$. On d\'efinit $$A=\begin{array}[t]{c}
\lim \vspace{-0.35 cm}\\ _{\longleftarrow} \end{array} \,(A_k,p_k)=\left\{(a_k)\in\prod_{k\ge 1} A_k,\ p_k(a_{k+1})=a_k\right\}\,.$$

Cet espace est compact, et on peut d\'efinir une distance ultra-m\'etrique $d_A$ comme suit. $$d_A((a_k),(b_k))=  \frac{1}{2^{\min\{j,\ a_j\ne b_j\}}}\,.$$

\begin{lemm}\label{lrep} Si $K$ est un compact de $Z$, alors il existe une application 
$I:A\to K$ surjective et $2$-Lipschitz.\end{lemm}

\Preuve Soit $\BBB^1=\{B_j^1\}_{1\le j\le N_1}$ un recouvrement par $N_1$ boules de rayon $1/2$. On d\'efinit $I_1(n)=x_n$, le centre de $B_n^1$.

On suppose que, jusqu'au rang $n\ge 1$, on a construit des applications $I_k:A_k\to K$ telles que

\bit

\item[(a)]  $I_k(A_k)$ est une $(1/2^k)$-approximation de $K$ (chaque point de $K$ est \`a distance
au plus $1/2^k$ de $I_k(A_k)$\,;

\item[(b)] Pour chaque $k$, $I_{k+1}(a)\in B(I_k p_k(a),1/2^k)$, pour tout $a\in A_{k+1}$.\eit

Pour chaque $a\in A_n$, $B_a=B(I_n(a),1/2^n)$ est recouvert par $N_{n+1}$ boules $$\BBB_a=\{B(x^a_j,1/2^{n+1})\}_{1\le j\le N_{n+1}}\,.$$

On d\'efinit $I_{n+1}(a,j)=x^a_j$. Notons que $x^a_j\in B_a$. 
Par d\'efinition, $I_{n+1}(\{a\}\times [1,N_{n+1}])$ est une $1/2^{n+1}$-approximation de $B_a$. Puisque
ces boules recouvrent $K$, on obtient (a) pour $k=n+1$. De m\^eme, (b) est v\'erifi\'e par construction puisque $p_n(a,j)=a$.

D'apr\`es (b), si $a=(a_k)\in A$, alors $(I_k(a_k))_k$ est une suite de Cauchy. 
On d\'efinit alors $$I(a)=\lim I_k(a_k)\,.$$

On note que $I(a)\in B(I_k(a_k),1/2^k)$ pour chaque $k$.

D'apr\`es (a) et (b), pour tout $x\in K$, on peut trouver une suite de boules dans chaque recouvrement qui contiennent $x$. Autrement dit,
$I$ est surjective.

Soient $a=(a_n)$ et $b=(b_n)$ des points de $ A$. On suppose que $d_A(a,b)= 1/2^{k}$. Du coup, $a_k=b_k$, et $I(a),I(b)\in B(I_k(a_k),1/2^k)$, donc
$|I(a)-I(b)|\le 1/2^{k-1}\le 2d_A(a,b)$. \endp

\bigskip

\begin{coro} Si $(X_n)$ est une suite de compacts de $Z$ qui tend vers un compact $X$, alors, quitte \`a extraire
une sous-suite, la suite $(I_n)$ d\'efinie par le Lemme \ref{lrep} tend vers une application $I:A\to Z$
telle que $I(A)=X$. \end{coro}

\Preuve Le lemme pr\'ec\'edent nous construit 
un espace m\'etrique $A$ et une suite d'applications $I^n:A\to Z$ $2$-lipschitziennes telles que $I^n(A)=X_n$. 
Cette suite est donc
\'equicontinue et le \th d'Ascoli nous extrait une sous-suite $(n_k)$ de sorte que $(I^{n_k})_k$
est convergente vers une application $2$-lipschitzienne 
$I:A\to Z$.
Il reste \`a v\'erifier que $X=I(A)$. Or, on a
 $$d_H(X_{n_k},I(A))\le \sup_{a\in A} |I^{n_k}(a)-I(a)|$$ donc $I(A)$ est la limite de Hausdorff de $(X_{n_k})$,
 ce qui \'etablit le corollaire, puisqu'il y a unicit\'e de la limite (parmi les ensembles ferm\'es).  \endp

\bigskip

{\noindent\bf Convergence des fonctions.} On \'etudie maintenant les limites de fonctions lorsque les espaces sont convergents.

\bigskip

On utilise la notion de module de continuit\'e pour traduire la notion de continuit\'e uniforme.

\begin{defi}[module de continuit\'e]
Si $f:X\to Y$ est continue, un module de continuit\'e est une fonction $\om:\R_+\to\R_+$ telle que $\om(t)\to 0$ avec $t$ et, pour tous $x,y\in X$,
$$|f(x)-f(y)|\le \om(|x-y|)\,.$$\end{defi}

On remarque aussi que si une famille de fonctions est uniform\'ement \'equicontinue, il existe un module de continuit\'e qui ne d\'epend pas de la fonction
de la famille.

\bigskip

\noindent{\sc Definition. ---} Soient $Z$ et $Y$ deux espaces m\'etriques compacts, 
$(X_n)$ une suite de
compacts de $Z$ qui converge vers un compact $X$. Soit $(f_n:X_n\to Y)$ une suite 
d'applications continues. 
On dit que $(f_n)$ converge vers $f:X\to Y$ si, pour tout
$\ep >0$ et pour tout $z\in X$, il existe $\al>0$, tel que, pour tout $n$ assez grand, 
si $z_n\in X_n$ v\'erifie
$d(z_n,z)<\al$ alors $d(f_n(z_n),f(z))<\ep$.

On obtient ainsi

\begin{theo}\label{Asc0}  Soient $Z$ et $Y$ deux espaces m\'etriques compacts, $(X_n)$ une suite de
compacts de $Z$ qui converge vers un compact $X$. Soit $(f_n:X_n\to Y)$ une suite d'applications continues 
qui admet un module de continuit\'e uniforme.

Il existe une sous-suite $(f_{n_k})$ qui converge vers une application continue $f:X\to Y$,
munie du m\^eme module de continuit\'e.

En particulier, si la suite $(f_n)$ est isom\'etrique, $L$-Lipschitz ou $\eta$-quasisym\'etrique, 
il en
est de m\^eme de la limite.\end{theo}

\Preuve Soit $\om$ le module de continuit\'e uniforme de $(f_n)$. 
On note $A$ le compact obtenu par le Lemme \ref{lrep}
pour  $X_n$ et on suppose que $I^n$ tend vers
$I$ uniform\'ement. 

On obtient ainsi une suite d'applications  $g_n=f_n\circ I^n:A\to Y$ uniform\'ement \'equicontinue.

Le \th d'Ascoli nous permet  d'extraire une  sous-suite convergente vers une application $g:A\to Y$. 

Soient $z,w\in X$. Il existe $a,b\in A$ tel que $I(a)=z$ et $I(b)=w$. Du coup,
$$|g(a)-g(b)|= |f_n(I^n(a))-f_n(I^n(b))|\pm (|g_n(a)-g(a)|+|g_n(b)-g(b)|)$$
donc $$|g(a)-g(b)|= \lim |f_n(I_n(a))-f_n(I_n(b))|\,.$$

Du coup, $g(a)=g(b)$ si $I(a)=I(b)$. Autrement dit, il existe
$f:X\to Y$ telle que $g=f\circ I$. Cette application $f$ a le m\^eme module de continuit\'e\,:
$$|f(z)-f(w)|= \lim |f_n(I^n(a))-f_n(I^n(b))|\le \lim \om(|I^n(a)-I^n(b)|)=\om(|z-w|)\,.$$
Le reste suit, et est laiss\'e au lecteur.
\endp

\end{document}